\documentclass[a4paper,11pt]{article}
\usepackage[vlined]{algorithm2e}

\SetFuncSty{textsc}
\SetKwFunction{frun}{Run}
\SetKwHangingKw{arun}{\frun}
\SetKwFunction{fset}{Set}
\SetKwHangingKw{aset}{\fset}
\SetKwFunction{fselect}{Select}
\SetKwHangingKw{aselect}{\fselect}
\SetKwFunction{fcompute}{Compute}
\SetKwHangingKw{acompute}{\fcompute}
\SetKwFunction{fsolve}{Solve}
\SetKwHangingKw{asolve}{\fsolve}
\SetKwFunction{festimate}{Estimate}
\SetKwHangingKw{aestimate}{\festimate}
\SetKwFunction{fmark}{Mark}
\SetKwHangingKw{amark}{\fmark}
\SetKwFunction{frefine}{Refine}
\SetKwHangingKw{arefine}{\frefine}
\SetKwFunction{compute}{Compute}
\SetKwFunction{set}{Set}

{\algorithm}%
{\endalgorithm}

\newcommand{\de}{\varepsilon}
\newcommand{\nvec}{\mathbf{n}}

\usepackage{amsmath,amsthm,amssymb,enumerate, amsbsy}
\usepackage{amssymb}
\usepackage[T1]{fontenc}
\usepackage{gensymb}
\usepackage{longtable}
\usepackage[square,sort&compress,comma,numbers]{natbib}
\usepackage[sc]{mathpazo}
\usepackage{multirow} 
\usepackage{newtxtext,newtxmath}
\usepackage{subfig}
\usepackage{graphicx}
\usepackage{epstopdf}
\usepackage[colorinlistoftodos]{todonotes}
\usepackage{epsfig, setspace, subfig}
\usepackage{hyperref}
\usepackage{cancel}
\usepackage[margin=2.5cm]{geometry}
\usepackage{tikz}
\usepackage[titletoc]{appendix}
\usepackage{caption}
\usetikzlibrary{shapes,calc}
\usepackage{verbatim}
\usepackage{mathrsfs}
\usepackage{accents}
\usepackage[utf8]{inputenc}
\usepackage{float}
\restylefloat{table}
\usepackage{multirow,diagbox,tabularx,blindtext,tabulary,tabularht,longtable,blkarray}

\usetikzlibrary{decorations.pathmorphing}
\usetikzlibrary{decorations.pathreplacing}
\usetikzlibrary{positioning}
\usetikzlibrary{shapes}
\usetikzlibrary{arrows}
\usetikzlibrary{patterns}
\usetikzlibrary{fadings}
\usetikzlibrary{plotmarks}
\usetikzlibrary{calc}
\usetikzlibrary{intersections}
\tikzstyle{every picture}+=[font=\footnotesize]
\usepackage{paralist}
\usepackage{empheq}
\usepackage{bbm}
\usepackage{latexsym}     
\usepackage{enumerate}
\usepackage{enumitem}
\usepackage{algorithm2e}

\setlist{noitemsep, topsep=0.8ex, partopsep=0pt
	, leftmargin=3em}
\setlist[1]{labelindent=\parindent}

\newlist{axioms}{enumerate}{1}
\setlist[axioms]{font=\bfseries}

\newlist{alphenum}{enumerate}{1}
\setlist[alphenum]{label=\textbf{(\alph*)}, leftmargin=4em}

\newlist{alphienum}{enumerate}{1}
\setlist[alphienum]{label=\textit{(\alph*)}}

\newlist{romanenum}{enumerate}{1}
\setlist[romanenum]{label=\textit{(\roman*)}}

\newlist{romaninenum}{enumerate*}{1}
\setlist[romaninenum]{label=\textit{(\roman*)}}
\usepackage[vlined]{algorithm2e}
\SetKwIF{If}{ElseIf}{Else}{if}{}{else if}{else}{endif}
\SetKwFor{For}{for}{}{endfor}
\usepackage[noabbrev, capitalise]{cleveref}

\theoremstyle{definition}

\theoremstyle{remark}

\numberwithin{equation}{section}

\newcommand{\Qvec}{\mathbf{Q}}
\newcommand{\abs}[1]{\vert #1 \vert}
\newcommand{\dx}{\,{\rm d}x}

\newcommand{\norm}[1]{{\vert\kern-0.25ex\vert #1 
		\vert\kern-0.25ex \vert}}

\title{Bayesian Parameter Identification in the Landau-de Gennes Theory for Nematic Liquid Crystals}
\author{\stepcounter{footnote} \stepcounter{footnote} \stepcounter{footnote} Heiko Gimperlein\footnote{ 
	Engineering Mathematics, 	University of Innsbruck, Innsbruck, Austria. Email:   heiko.gimperlein@uibk.ac.at} $\;\;$
    Ruma R. Maity\footnote{ 
	Engineering Mathematics, 	University of Innsbruck, Innsbruck, Austria. Email:  ruma.maity@uibk.ac.at} $\;\;$
	Apala Majumdar \footnote{ Department of Mathematics and Statistics, University of Strathclyde, Glasgow, United Kingdom. Email: apala.majumdar@strath.ac.uk}        
	$\;\;$
    Michael Oberguggenberger\footnote{Engineering Mathematics, 	University of Innsbruck, Innsbruck, Austria. Email. michael.oberguggenberger@uibk.ac.at}
	\date{}
}

\begin{document}

\maketitle
\begin{abstract} 
This manuscript establishes a pathway to reconstruct material parameters from measurements within the Landau-de Gennes model for nematic liquid crystals. We present a Bayesian approach to this inverse problem and analyse its properties using given, simulated data for benchmark problems of a planar bistable nematic device. In particular, we discuss the accuracy of the Markov chain Monte Carlo approximations, confidence intervals and the limits of identifiability.
\end{abstract}

\section{Introduction}

Nematic liquid crystals (NLCs) are classical mesophases that are more ordered than conventional liquids and less ordered than conventional solids \cite{dg}. The constituent asymmetric NLC molecules move freely but tend to align along locally preferred directions, referred to as nematic directors. Consequently, NLCs exhibit direction-dependent responses to external stimuli, e.g., light, electric fields, temperature and this directionality and material softness make them excellent candidate working materials for electro-optic devices, photonics, elastomers and biomimetic materials \cite{lagerwallreview}.

NLC applications crucially depend on a comprehensive understanding of the following question: given a specified NLC system, can we mathematically predict and control the experimentally observable NLC configura\-tions? There are multiple mathematical theories for NLCs, that describe both the equilibrium and non-equilibrium  
 properties of NLCs in confinement, across different length and time-scales. For example, there are molecular models that incorporate information about molecular shapes, sizes and intermolecular interactions; there are mean-field models that average the intermolecular interactions in terms of a mean field acting on the molecules and retain some molecular-level information; there are macroscopic continuum theories which describe the NLC state in terms of a macroscopic order parameter with no explicit connections to the underlying molecular-level or microscopic details \cite{ballmajumdar2010}. We focus on the powerful continuum Landau-de Gennes (LdG) theory for NLCs in this paper, which is a phenomenological variational theory that has enjoyed tremendous success in the modelling and applications of NLCs across disciplines \cite{majumdarEJAM2010}. The LdG theory describes the NLC state in terms of a macroscopic order parameter, the $\Qvec$-tensor order parameter which contains information on the nematic directors and the degree of orientational ordering about the directors \cite{majumdarEJAM2010}. As with many variational theories in  materials science, there is a LdG free energy which is a measure of the free energy of a NLC system and the modelling hypothesis is that physically observable equilibrium NLC configurations are mathematically modelled by LdG energy minimisers, subject to the imposed boundary conditions. The energy minimisers, and indeed, all critical points of the LdG free energy are classical solutions of the associated Euler-Lagrange equations, which are systems of nonlinear partial differential equations (PDE) with (typically) multiple solutions \cite{majumdarzarnescu2010}. Therefore, the equilibrium properties of NLC systems or NLC equilibria are studied in terms of boundary-value problems for the Euler-Lagrange equations of the LdG free energy, which is a complex PDE problem in its own right. The LdG free energy, and consequently, the Euler-Lagrange equations have multiple phenomenological parameters related to the temperature, material properties, elastic constants, chirality etc. There are inherent uncertainties in these phenomenological parameters, the experimental system, in the mathematical model and the numerical methods. We need uncertainty quantification (UQ) to (i) estimate the uncertainties in the mathematical predictions and control the errors and (ii) understand the limitations of the models themselves. One can approach UQ in at least two different ways: forward UQ which focuses on how uncertainties in the LdG model inputs propagate to the model outputs, e.g., the predicted values of the equilibrium LdG $\Qvec$-tensor order parameter so that we can compute error bounds for the LdG model outputs and assess their accuracy/reliability. Inverse UQ focuses on the question -- given a LdG Q-tensor or a distribution for the same, perhaps constructed from experimental data, can we estimate the unknown/uncertain LdG model inputs for benchmark problems? 
 
 We focus on inverse UQ problems in the LdG theory in this manuscript. Our goal is to establish an algorithmic pathway from optical measurements and experimentally measured dielectric data to experimental measurements of the $\Qvec$-tensor order parameter (this is already known in the literature), interpret the experimental measurements of the $\Qvec$-tensor order parameter as an equilibrium measurement of the LdG $\Qvec$-tensor order parameter, and then use inverse UQ to reconstruct the LdG model inputs from the observed values/given measurements of the LdG $\Qvec$-tensor order parameter. We focus on solving this inverse problem using statistical tools --  Bayes' theorem \cite{Stuart_2010_acta_numerica, DashtiStuart:2017Bayesian, HeldBove:2014} and Markov chain Monte Carlo sampling \cite{Roy2020, Roberts1998, Gamerman2006}. Our contributions are summarized as follows.
\begin{itemize}
    \item A comprehensive formulation of the inverse problems that can estimate the LdG model inputs from appropriate NLC experimental measurements.
    \item  Bayesian inference to reconstruct LdG model inputs from given data on the observed values of the LdG $\Qvec$-tensor order parameter, using Markov chain Monte Carlo (MCMC) methods. 
    \item  Computation of posterior distributions of the LdG model inputs, given data and distributions for the LdG $\Qvec$-tensor order parameters for a benchmark example, accompanied by computations of Bayesian estimators (posterior mean, posterior median and corresponding  confidence intervals for the parameters of interest).
\end{itemize}
There are existing well established methods for measuring the phenomenological parameters in the LdG free energy or the LdG model inputs, e.g., one can use the nematic-isotropic transition temperature to estimate some bulk LdG constants; researchers use the classical Freedericksz transition or classical NLC instabilities, on exposure to external electric fields to measure NLC elastic constants -- see \cite{Klus:14} wherein this method has been applied to measure the elastic constants of the common NLC materials, 6CHBT and E7;  a dielectric and an optical method has been used in \cite{Faetti2009} to measure the splay and bend elastic constants  for the  Merck nematic mixture E49;  an optical method is used in \cite{Bankova2024} for measuring
the twist elastic constant along with the splay and bend
elastic constants of  commonly used NLCs: 5CB, 6CHBT, and E7.
Recently a neural networks-based method has been explored in \cite{Zaplotnik2023} to determine NLC elastic constants based on combined modelling
of  NLC dynamics,  light transmission, and  supervised machine learning. Other relevant works on Bayesian inference for inverse problems include identification/reconstruction of the optical parameters in  the Ginzburg-Landau equation,  which is  a  celebrated fundamental model for dissipative optical solitons \cite{Zhao2024}. 

We focus on a benchmark example motivated by the planar bistable nematic device reported in \cite{Tsakonas}. The experimental domain is a shallow three-dimensional square well with tangent boundary conditions on the well surfaces, that constrain the NLC molecules to lie on the well surfaces.  This device has been experimentally manufactured and is known to support multiple competing NLC equilibria (referred to as diagonal and rotated solutions), without any external electric fields, if the square well is sufficiently large in size. This device has been modelled in \cite{Tsakonas} and \cite{MultistabilityApalachong}, in the LdG framework. In both cases, the authors neglect the height of the square well (since the height is typically much smaller than the cross-sectional dimensions) and take the computational domain to be a square domain, subject to tangent boundary conditions for the LdG $\Qvec$-tensor order parameter on the square edges. They work in a reduced two-dimensional LdG framework (to be described in Subsection \ref{Inverse-problem-formulation}\ref{sec: Q to para}), so that the governing Euler-Lagrange equations are a system of two nonlinear PDEs, with two model parameters (denoted by $\alpha$ and $\beta$). The model parameters contain information about the square edge length, the temperature and the NLC material properties. The authors of \cite{Tsakonas} and \cite{MultistabilityApalachong} work in a deterministic framework, assume $\alpha$ and $\beta$ are known, and recover the diagonal and rotated solutions for appropriate choices of $\alpha$ and $\beta$ in the reduced LdG framework. Of course, they also study additional aspects such as how multistability (co-existence of multiple stable diagonal and rotated solutions) depends on the square edge length and the boundary conditions, but there is no uncertainty in their work. In \cite{dalby2024}, the authors study the solution landscape of the LdG Euler-Lagrange equations, for this benchmark problem, with additive and multiplicative noise whereby this noise captures uncertainties in a holistic sense. The authors conclude that the deterministic LdG predictions are fairly robust, notwithstanding the addition of stochastic noise terms, for physically relevant values of $\alpha$ and $\beta$.

In real life, $\alpha$ and $\beta$ are almost certainly not known exactly. In this paper, we focus on reconstructing $\alpha$ and $\beta$, given diagonal and rotated solutions for this benchmark problem on a square domain or equivalently, given $\Qvec$-tensor solutions of the reduced LdG Euler-Lagrange equations for this benchmark problem. We treat the diagonal and rotated solutions as given observed solutions, $\Qvec_{\rm obs}$, which are exactly known. We also assume that the square edge length and the temperature are exactly known. We assume prior distributions for $\alpha$ and $\beta$, for which there is no general consensus. Hence, we experiment with two different prior distributions for $\alpha$ and $\beta$ -- a uniform or diffusive prior and a suitably constructed Gaussian prior distribution. We use Bayesian methodology to construct the posterior distribution for $\alpha$ and $\beta$ (the probability distribution for $\alpha$ and $\beta$, given $\Qvec_{\rm obs}$), given a prior distribution for the LdG model inputs and a suitably constructed likelihood function (the probability distribution of $\Qvec_{\rm obs}$, given the LdG model inputs). The posterior distributions are computed using MCMC methods, implemented by a Metropolis-Hastings algorithm. Our methods work for several physically relevant values of $\alpha$ and $\beta$, as illustrated by computations of Bayesian estimators and confidence intervals for the estimators. We also identify situations for which our methods do not work and offer some heuristic insights into the possible reasons for non-identifiability.

We believe that our work is a good attempt at incorporating UQ into the LdG theoretical framework, which can much improve modelling capabilities. The liquid crystal research community needs multiple reliable and efficient methods for estimating material properties, for accurate and powerful mathematical modelling. For example, one could partially use our methods to compute $\alpha$ and $\beta$ (with error estimators) from a set of experimentally derived NLC director patterns, for a given NLC material. Once $\alpha$ and $\beta$ are reliably known for a NLC material, researchers could numerically canvass complex solution landscapes of the LdG Euler-Lagrange equations on complex geometries and accurately predict equilibrium NLC properties or reciprocally, design new NLC systems that can yield desired equilibrium properties. We do not argue that inverse UQ will overtake conventional methods for estimating NLC material properties, but it can provide a decent alternative and inverse UQ also has great potential for machine learning approaches to liquid crystal research.

This paper is organized as follows. In Section \ref{Inverse-problem-formulation}, we outline the pathway from optical and dielectric measurements to estimates of the LdG model inputs, in a reduced two-dimensional setting. In Section~\ref{sec:Bayes}, we outline our Bayesian methodology. In Section~\ref{sec:computational_framework}, we describe our  finite element framework for the  benchmark example motivated by the planar bistable nematic device reported in \cite{Tsakonas}. In Sections~\ref{numerical-constant-L} and \ref{numerical-two-para}, we present several numerical experiments on the reconstruction of posterior distributions of $\alpha$ and $\beta$, given $\Qvec_{\rm obs}$ and information about temperature and square edge length. In Section~\ref{Disc}, we discuss non-identifiability or situations wherein our methods do not converge, which might necessitate further study. We conclude with some perspectives in Section~\ref{conclusions}.

\section{An inverse problem in the Landau-de Gennes framework}\label{Inverse-problem-formulation}
In this section, we outline the main steps involved in determining the nematic liquid crystal (NLC) material parameters from experimental measurements, e.g, optical measurements. The key experimental quantity is the dielectric anisotropy, $\varepsilon$, which is a measure of the NLC response to external electric fields. The primary mathematical variable is the Landau-de Gennes (LdG) $\Qvec$-tensor order parameter, which is a symmetric traceless $3 \times 3$ matrix that encodes information about the nematic directors and the order parameters. The LdG $\Qvec$-tensor order parameter can be empirically related to $\varepsilon$, providing a link between experiments and theory. Equally, the LdG $\Qvec$-tensor order parameter depends phenomenologically on various NLC material parameters and the temperature, as dictated by the LdG modelling framework. Hence, we can, in principle, compute the NLC material parameters from experimental measurements of $\varepsilon$ in the following three steps:
\begin{enumerate}
    \item inferring the dielectric parameters embedded in $\varepsilon$ from optical measurements;
    \item determining the LdG $\Qvec$-tensor from $\varepsilon$;
    \item and using the LdG framework to determine the NLC material parameters from measurements/computations of the equilibrium values of the LdG $\Qvec$-tensor order parameter.
\end{enumerate}
This paper is primarily devoted to the last step in the above sequence, but we briefly describe the first two steps in the remainder of this section.

\begin{figure}
    \centering
    \includegraphics[width=0.9\linewidth]{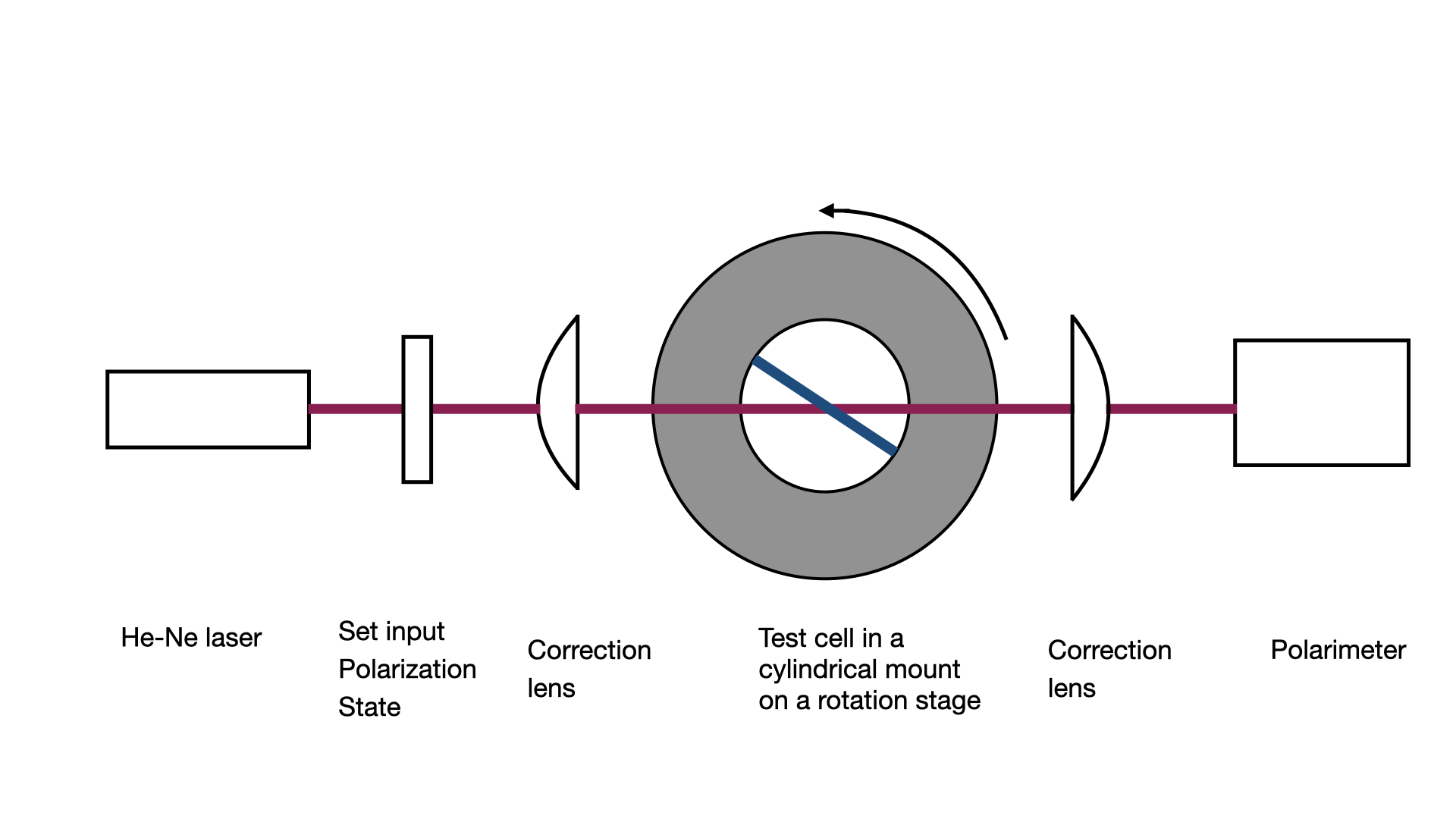}
    \caption{Experimental setup used at Hewlett Packard Laboratories, Bristol, as reported in \cite{ln07}.}

    \label{fig:measurement-set-up}
\end{figure}
\subsection{From optical measurements to $\varepsilon$ }
\label{sec:optical}

Lionheart and Newton \cite{ln07} describe a typical experimental set-up in three dimensions (3D). As illustrated in Figure \ref{fig:measurement-set-up}, an NLC sample is illuminated by a monochromatic, polarized laser beam. The  normalized Stokes parameters of the transmitted light are measured depending on the incident angle and polarization. The measured data is then compared to the results of a forward model for the light transmission in the sample. The key model parameter is the dielectric anisotropy tensor, $\varepsilon$, which is fitted to the data. We briefly review the forward model for light transmission from Section 3 in \cite{ln07}. 

We assume that the direction of the laser beam is along the $z$-axis, and its polarization defines a vector in the $xy$-plane, perpendicular to the direction of the beam. The monochromatic electric field $\mathbf{\mathcal{E}}$ \cite{Robson:1974} is of the form $\mathbf{\mathcal{E}}=(\mathcal{E}_x, \mathcal{E}_y)$ with $\mathcal{E}_x:= a \cos(\omega t -kz)$ and  $\mathcal{E}_y:= b \cos(\omega t -kz+\delta)$. Here $\omega$ is the angular frequency, $k$ is the wave number, and $\delta$ is the phase difference between the two components. The experimental measurements yield values for the  Stokes parameters $S=(S_0,S_1,S_2, S_3)$ of the electric field, defined by 
\begin{align*}
    S_0:=a^2+b^2, \quad S_1:=a^2-b^2, \quad S_2:=2ab \cos\delta, \quad S_3:=2ab \sin\delta.
\end{align*}
These parameters characterize the state of the polarization, which solely depends on the amplitudes $a$, $b$, and the phase difference $\delta$. 
The next step is to consider the interaction between the NLC and propagating electric and magnetic fields, $\mathbf{\mathcal{E}}$ and $\mathbf{\mathcal{H}}$. In a NLC medium, the Berreman  field vector $X:=(\mathcal{E}_x,\mathcal{H}_y, \mathcal{E}_y,-\mathcal{H}_x )^T$ evolves according to the linear differential equation 
\begin{align}\label{Berreman-equn}
    \frac{\partial}{\partial z}X= -\frac{i \omega }{c}M {X}.
\end{align}
Here, the Berreman matrix $M$ depends on the dielectric tensor $\varepsilon = (\varepsilon_{ij})$ as shown below \cite{Berreman:72}: \begin{equation}\label{BerremanM}M =
\begin{pmatrix} -\frac{\varepsilon_{13}}{\varepsilon_{33}} \xi&\mu_0 c\frac{\varepsilon_{33}-\xi^2}{\varepsilon_{33}}&- \frac{\varepsilon_{23}}{\varepsilon_{33}}\xi&0\\
\varepsilon_0 c\left(\varepsilon_{11}-\frac{\varepsilon_{13}^2}{\varepsilon_{33}}\right)&-\frac{\varepsilon_{13}}{\varepsilon_{33}} \xi&\varepsilon_0 c\left(\varepsilon_{12}-\frac{\varepsilon_{13}\varepsilon_{23}}{\varepsilon_{33}}\right)&0\\
0&0&0&\mu_0c\\
\varepsilon_0 c\left(\varepsilon_{12}-\frac{\varepsilon_{13}\varepsilon_{23}}{\varepsilon_{33}}\right)&-\frac{\varepsilon_{23}}{\varepsilon_{33}} \xi&\varepsilon_0 c\left(\varepsilon_{22}-\frac{\varepsilon_{23}^2}{\varepsilon_{33}}-\xi^2\right)&0
\end{pmatrix}.\end{equation}
The permeability of the surrounding medium (typically air) is denoted by $\mu_0$, $c$ is the speed of light, and $\xi$ is determined by the incident angle and the refractive index of the surrounding medium. Moreover, one can show that $\varepsilon$ is uniquely determined by the matrix $M$.

The forward problem can thus be written as $S = F(\varepsilon,\mu_0,c,\xi)$, with $S$ the vector of Stokes parameters. The value of $S$ can be computed from the forward model, as a function of $(\varepsilon,\mu_0,c,\xi)$, or measured experimentally.
The inverse problem is to compute $\varepsilon$ from measurements of $S$ for various $\xi$, given $\mu_0$, $c$, as studied in \cite{ln07}. In the remainder of the manuscript, we assume that $\varepsilon$ has been determined from measurements of $S$.

\subsection{From $\varepsilon$ to the LdG $\Qvec$-tensor order parameter}
\label{sec:de}

The dielectric anisotropy tensor, $\de$, is a measure of the anisotropic NLC response to external electric fields \cite{dg}. The tensor, $\de$, has two key components: $\de_{||}$ which measures the dielectric response along the nematic (NLC) director (the distinguished direction of averaged molecular alignment), and $\de_{\perp}$ which measures the dielectric response in all orthogonal directions to the director \cite{dg}. If $\de_{||} - \de_{\perp} >0$, then the NLC molecules reorient to align with the applied external electric field and if $\de_{||} - \de_{\perp} <0$, then the NLC molecules reorient to be orthogonal to the applied external electric field.

The dielectric anisotropy tensor, $\de$, can be measured with laboratory studies \cite{patranabish2021} (as described in Subsection~\ref{sec:optical}) and is often used as an implicit phenomenological definition for the Landau-de Gennes (LdG) $\Qvec$-tensor order parameter, i.e., $\de$ and $\Qvec$ are related to each other as shown below  \cite{calderer}:
\begin{equation}
    \label{eq:de}
\varepsilon_{ij} = \frac{1}{3}\left(\varepsilon_{\parallel} + 2 \varepsilon_{\perp}\right) \delta_{ij} + (\varepsilon_{\parallel} - \varepsilon_{\perp})
Q_{ij} = \frac{1}{3}\ \mathrm{tr}\, (\varepsilon) \ \delta_{ij} + (\varepsilon_{\parallel} - \varepsilon_{\perp})
Q_{ij}; ~ i,j=1\ldots 3,
\end{equation}
where $\varepsilon_{\perp}$ and $\varepsilon_{\parallel}$ are the dielectric permittivities perpendicular and parallel to the director respectively, measured in units of the vacuum permittivity $\varepsilon_{0}$. In (\ref{eq:de}), there is an implicit assumption that the NLC sample under consideration is uniaxial, i.e., the NLC sample has a single distinguished director so that all directions perpendicular to the NLC director are physically equivalent.

Next, we describe a thought experiment to measure the dielectric permittivities in the matrix equation above. One could take a thin slab-based 3D geometry, wherein the top and bottom surfaces are uniformly rubbed to fix the NLC director in the same direction, on the top and bottom surfaces. In the absence of any external constraints, we expect the NLC director to be uniform throughout the domain, as dictated by the boundary treatments. The electric field can be applied in the direction of the boundary condition, or equivalently in the direction of the nematic director for spatially homogeneous uniform samples. One can experimentally measure the dielectric responses parallel and perpendicular to the applied electric field, to measure $\de_{||}$ and $\de_{\perp}$ respectively. Once there are measurements of $\de_{||}$ and $\de_{\perp}$ in (\ref{eq:de}) for a given NLC material, we can experimentally measure the dielectric anisotropy tensor for generic confined NLC systems (using the methodology outlined in Subsection~\ref{sec:optical}) and use (\ref{eq:de}) to reconstruct the corresponding $\Qvec$. To summarise, (\ref{eq:de}) is a matrix equation relating the nine matrix components of the dielectric tensor, $\de$, and the LdG $\Qvec$-tensor, where $\de$ can be experimentally measured and $\de_{||}$ and $\de_{\perp}$ are given from prototype experiments on spatially homogeneous NLC systems.

\subsection{From the LdG $\Qvec$-tensor order parameter to the NLC material parameters}\label{sec: Q to para}
Our main goal is to provide algorithms that can compute the LdG model parameters from given data for the LdG $\Qvec$-tensor order parameters.

The LdG theory is the most powerful continuum theory for NLCs in the literature \cite{dg}. The LdG theory describes the NLC state by a macroscopic order parameter, the LdG $\Qvec$-tensor order parameter which contains information about the NLC directors and the degree of orientational ordering about them \cite{Mottram2014IntroductionTQ}. The LdG $\Qvec$-tensor is a symmetric traceless $3 \times 3$ matrix whose eigenvectors model the NLC directors and the corresponding eigenvalues measure the degree of orientational ordering about the directors. More precisely, the admissible LdG $\Qvec$-tensors belong to the space $S_0^n:=\{ \mathbf{P}= (P_{ij})_{n \times n}|\, P_{ij}=P_{ji}, \textrm{tr} \,\mathbf{P}=0\}$ for $n=2,3$.

The NLC phases can be categorised as follows, depending on the eigenvalue structure of the $\Qvec$-tensor: (i)  $\Qvec =0$ for isotropic (disordered) NLC phases with no defined special material directions; (ii) uniaxial if $\Qvec$ has two degenerate non-zero eigenvalues and the NLC director is the eigenvector corresponding to the non-degenerate eigenvalue and (iii) biaxial if $\Qvec$ has three distinct eigenvalues $\lambda_1> \lambda_2 > \lambda_3$, and there are two NLC directors corresponding to the two largest eigenvalues. For uniaxial NLC phases, the $\Qvec$-tensor can be written as
\[
\Qvec:=s ( \mathbf{n}\otimes  \mathbf{n} -\mathbf{I}/3 )
\]
where $\mathbf{I}$ is the $3 \times 3$ identity matrix  and  $\mathbf{n}$ is the leading eigenvector, with the non-degenerate eigenvalue, that models the single distinguished direction of averaged molecular alignment
at every point in space. The scalar order parameter, $s$, measures the degree of orientational order about $\mathbf{n}.$ The defect set is often identified with the nodal set of $s$, wherein the NLC ordering breaks down.

\medskip

In the absence of surface energies and external fields, we work with the simplest form of the 3D LdG energy \cite{majumdarzarnescu2010}:
\begin{equation}
\label{eq:energy}
    \mathcal{E}(\Qvec):=\int \Big( \frac{L}{2}|\nabla \Qvec|^2 + f_B(\Qvec) \Big) \, \dx, \ \ \text{where }\ f_B(\Qvec):=  \frac{A}{2} \textrm{tr} \, \Qvec^2 - \frac{B}{3}\textrm{tr} \,\Qvec^3 +\frac{C}{4} (\textrm{tr} \,\Qvec^2)^2 .
\end{equation}

Here, $L>0$ is a material-dependent elastic constant, the Dirichlet elastic energy density penalises spatial inhomogeneities and $f_B(\Qvec)$ is the bulk potential \cite{majumdarEJAM2010}.
The variable $A:=\alpha_0 (T -T^*)$ is a rescaled temperature where $\alpha_0, B,C>0$ are material dependent constants \cite{dg}, and $T^*$ is the characteristic nematic supercooling temperature. The critical points and minimisers of the bulk potential can be explicitly computed in terms of an algebraic problem \cite{majumdarEJAM2010}. The critical points of $f_B$ are either uniaxial $\Qvec$-tensors or the isotropic phase, $\Qvec=0$. The rescaled temperature $A$ has three characteristic values: 
\vspace{-0.4cm}
\begin{enumerate} \item $A=0$, below which the
isotropic phase $\Qvec =0$ loses stability; \item the nematic-isotropic transition temperature, $A=\frac{B^2}{27C}$, at which $f_B$ is minimised by the isotropic phase and a continuum of uniaxial states with $s=s_+=\frac{B +\sqrt{B^2 - 24 AC}}{4C}$ and $\mathbf{n}  $ arbitrary;
and \item  the nematic super-heating temperature,  $A=\frac{B^2}{24C}$ above which the isotropic state is the unique critical point of $f_B.$
\end{enumerate}
\vspace{-0.2cm}

There are some reported  values  of the bulk constants \cite{1975Itlc} for the canonical NLC material, MBBA, which are $B=0.64\times 10^6 \, J/m^3$ and  $C=0.35 \times 10^6 \, J/m^3$, $\alpha_0=0.042 \times 10^6 \, J m^{-3} K^{-1}$, $T^*= 45^o C, T_c=46^o C$, where $T_c$ is the nematic-isotropic transition temperature. Typical values of $L$ are around $10^{-12} - 10^{-11} J/m$ \cite{Karlj_Majumdar2014}. The physically observable configurations are modelled by minimisers of the LdG energy (\ref{eq:energy}), subject to the imposed boundary conditions. The energy minimisers, and indeed all critical points of the LdG energy in (\ref{eq:energy}) are analytic solutions of the corresponding Euler-Lagrange equations, which are a system of five nonlinear and coupled elliptic partial differential equations \cite{majumdarzarnescu2010}:
\begin{equation}
\label{eq:EL}
L \Delta \Qvec =  A \Qvec - B \left(\Qvec \Qvec - \frac{1}{3} \left(\textrm{tr} \,\Qvec^2 \right) \mathbf{I} \right) + C \left(\textrm{tr} \,\Qvec^2\right) \Qvec.
\end{equation}

For a given non-homogeneous Dirichlet boundary condition $ \Qvec_b\in  {H}^\frac{1}{2}(\partial\Omega;S^3_0),$ and the admissible set,  $\mathcal{A}(\Qvec_b ):=\{ \Qvec \in {H}^1(\Omega; S_0^3)|\,\, \Qvec =  \Qvec_b \text{ on }{\partial\Omega}\}$,   there exists an energy minimiser, i.e., the optimisation problem
\[ \min_{\Qvec \in \mathcal{A}(\Qvec_b )} \mathcal{E}[\Qvec] \] has a solution from the direct methods in the calculus of variations \cite{Gartland1998}.
 Indeed, the bulk density $f_B(\bullet)$ satisfies the growth condition \cite[Corollary 4.4]{Gartland1998}, $ f_B(\Qvec) \geq C +C |\Qvec|^2$ for all $\Qvec\in S_0^3 $ for some constants $C >0.$
The coercivity and convexity of the elastic energy in $\nabla \Qvec,$
together with the above mentioned growth condition of $f_B$, suffice to guarantee the existence of a global LdG energy minimiser in $\mathcal{A}(\Qvec_b )$.

We work in a two-dimensional (2D) framework in this paper. Using the gamma convergence arguments in \cite{golovatymontero2015}, one can argue that for thin-film geometries (where the height is much smaller than the cross-sectional dimensions) and for certain types of physically relevant boundary conditions, the energy minimisers of (\ref{eq:energy}) have a fixed eigenvector with an associated fixed eigenvalue (determined by the temperature and the imposed boundary conditions). In other words, one can define a reduced LdG $\Qvec$-tensor in 2D scenarios, where
\begin{equation}
    \Qvec=s\left(2\nvec\otimes\nvec-\mathbf{I}_2\right).\label{eq:Q}
\end{equation}
Here, $\Qvec$ is independent of $z$ or the spatial coordinate along the height of the well, $\nvec\in\mathbb{S}^1$ is the nematic director in the $xy$-plane, $s$ is the scalar order parameter which measures the degree of nematic ordering about $\nvec$, and $\mathbf{I}_2$ is the $2\times 2$ identity matrix. The 2D nematic director can be interpreted as the preferred direction of averaged molecular alignment in the $xy$-plane, and the defect set is simply the nodal set of $s$, denoted by $\mathcal{S}$.
The symmetry and tracelessness of $\Qvec$ imply that there are only two independent components, $Q_{11}$ and $Q_{12}$, given by
        \begin{equation*}
            Q_{11} = s\cos 2 \vartheta,\quad Q_{12} = s \sin 2\vartheta,
        \end{equation*}
        where $\nvec = \left(\cos \vartheta, \sin \vartheta \right)$ and $\vartheta$ denotes the angle between $\nvec$ and the $x$-axis. The reduced LdG tensor can be related to the full 3D LdG $\Qvec$-tensor as follows:
\begin{equation}\label{eq:3D2D}
    \Qvec_f=
            \begin{pmatrix}
            Q_{11}-q_3 & Q_{12} & 0\\
            Q_{12} & -Q_{11}-q_3 & 0\\
            0 & 0 & 2q_3
        \end{pmatrix},
\end{equation} where $Q_{11}$ and $Q_{12}$ are as above and $q_3$ is a fixed known constant determined by the temperature and boundary conditions (coded in terms of surface energies). Note that $\Qvec_f$ is not, in general, an exact solution of (\ref{eq:EL}) but a good approximation to stable energy minimising solutions of (\ref{eq:EL}) in the thin-film limit.

The correspondence in \eqref{eq:3D2D} is only valid for 2D scenarios. In \cite{Han_majumdar_Zhang_2020}, the authors show that for the special temperature $A=-\frac{B^2}{3C}<0$, there exists a branch of exact 3D solutions of (\ref{eq:EL}) of the form (\ref{eq:3D2D}), so that the results in the 2D framework are  transferable to fully 3D scenarios. Hence, in the remainder of this manuscript, we focus on this special temperature for which $s_+ = \frac{B}{C}$ and work in the reduced LdG framework, with the reduced LdG tensor defined in (\ref{eq:Q}). The full 3D $\Qvec$-tensor, $\Qvec_f$, can be reconstructed from the reduced LdG-tensor as in (\ref{eq:3D2D}). The corresponding reduced 2D LdG free energy for $\Qvec:=(Q_{11},Q_{12})$, at the temperature $A =-\frac{B^2}{3C}$, is given by
\begin{equation}\label{eq:2Denergy}
  \mathcal{E}(\Qvec):=\int_\Omega\Big(\frac{L}{2}|\nabla \Qvec|^2  -\frac{B^2}{4C}\textrm{tr} \,\Qvec^2  +\frac{C}{4}(\textrm{tr} \,\Qvec^2)^2 \Big) \dx,
\end{equation}
where $\Omega$ is the 2D Lipschitz domain and the other quantities, $L, B$ and $C$ have been defined above. Let $\lambda$ be a characteristic length associated with $\Omega$; we can rescale the system using the change of variable $\tilde{x} = \frac{x}{\lambda}$ and drop all tildes in the subsequent rescaled/dimensionless formulation.
The associated rescaled Euler-Lagrange equations are \cite{Han_majumdar_Zhang_2020}
\begin{align}\label{eqn:constant-L}
\begin{split}
       \Delta Q_{11} = \frac{2C \lambda^2}{L} \Big(  Q_{11}^2 +Q_{12}^2 -\frac{B^2}{4C^2}\Big) Q_{11} ,\\
   \Delta Q_{12} = \frac{2C \lambda^2}{L} \Big(  Q_{11}^2 +Q_{12}^2 -\frac{B^2}{4C^2}\Big) Q_{12}.
\end{split}
\end{align}
Using the aforementioned values of the parameters $B,C$ for MBBA, we estimate $\frac{B^2}{4C^2}=0.83592.$

\noindent By rearranging the parameters in \eqref{eqn:constant-L}, we obtain the equivalent system of two coupled and nonlinear partial differential equations:
\begin{align}\label{pde-two-parameters}
    \begin{split}
         -\alpha \Delta \Qvec +   (\abs{\Qvec}^2-\beta)\Qvec = \,& 0 \text{ in }\Omega,\\
    \Qvec=\, &\Qvec_b \text{ on }\partial \Omega,
    \end{split}
\end{align}
where $\alpha:= \frac{L}{2C\lambda^2}$ and $\beta:=\frac{B^2}{4C^2}.$
If the temperature,  $A=-\frac{B^2}{3C}$, and the domain length-scale $\lambda$ are given, the parameters $\alpha,\beta$ determine $C=-\frac{3A}{4\beta}$, which allows us to compute $B=\sqrt{4C^2 \beta}$  and $L=2\alpha C \lambda^2$. 

Therefore, for a given temperature and domain length, it is sufficient to identify the parameters $\alpha$ and $\beta$ to determine the material parameters $L, B,C.$ For 2D polygons, $\lambda$ can be the edge length, usually in the range $\lambda \in \left[10^{-8}, 10^{-6}\right]m$, $C \in [10^5, 10^6] J m^{-3}$ and $L \in [10^{-12}, 10^{-11}] J m^{-1}$. Hence, it is reasonable to work with $\alpha \in \left(10^{-4}, 10^{-2} \right)$ as a representative range, along with $\beta \in \left(0.1, 1 \right)$. In what follows, we reconstruct $\alpha$ and $\beta$ in these physically relevant ranges from given solutions of the reduced LdG Euler-Lagrange equations (\ref{pde-two-parameters}) on a square domain, subject to Dirichlet boundary conditions for $Q_{11}$ and $Q_{12}$, for $A=-\frac{B^2}{3C}$.

If we work with arbitrary low temperatures $A<0$ (assumed to be known) in a strictly 2D scenario, then we do not have enough information to compute $B$. In such cases, $\alpha = \frac{L}{2\lambda^2 C}$ and $\beta = \frac{|A|}{2C}$ and one can reconstruct $C$ and $L$ from $\alpha$ and $\beta$, for given values of $A$ and $\lambda$.

\section{Bayesian methodology}
\label{sec:Bayes}

\subsection{Statistical inverse problems}
\label{subsec:inverse}

Consider the forward problem defined by an input-output map $\mathcal{F}$,
$$Y = \mathcal{F}(X,\Theta).$$
Here $X$ is the model input, $\Theta$ are model parameters and $Y$ is the model output. In the LdG context, $X$ is the information about characteristic domain/geometric length scales and ambient temperature (both assumed to be given), $\Theta$ are the parameters $\alpha$ or $(\alpha,\beta)$,
$Y$ is the bivariate output $\Qvec = (Q_{11},Q_{12})$, evaluated at the grid points of a shape  regular triangulation, and $\mathcal{F}$ is the map
$(\alpha, \beta)\to \Qvec$ defined through equation \eqref{pde-two-parameters}.
The inverse problem is to determine an estimate of the parameters, $\widehat{\Theta}$ that minimises the distance of the corresponding model output $\mathcal{F}(X, \widehat{\Theta})$ to the observed data $Y_{\rm obs}$, for given model input $X$.

We use a Bayesian approach to this inverse problem (as described in detail in \cite[Chapter 3]{Kaipio:2005}): $\Theta$ and $Y$ are viewed as random variables,  $X$ is known and the solution of the inverse problem is to compute the (posterior) probability distribution of $\Theta$. The Bayesian set-up requires (a) the so-called prior distribution $\pi_{\text{prior}}(\Theta)$ which encodes theoretical knowledge about the distribution of $\Theta$, (b) observed data $Y_{\rm obs}$ and (c) the likelihood function
$\pi(Y_{\rm obs}|\Theta)$, that is the probability distribution of observing $Y_{\rm obs}$ when $\Theta$ is given. The so-called posterior distribution $\pi_{\text{post}}(\Theta|Y_{\text{obs}})$ is the sought-after distribution of $\Theta$, given the observed
data $Y_{\text{obs}}$. The posterior distribution is computed by means of Bayes' theorem as shown below:
\begin{align}\label{posterior-distribution}
    \pi_{\text{post}}(\Theta):= \pi(\Theta| Y_{\text{obs}})=\frac{\pi_{\text{prior}}(\Theta) \pi( Y_{\text{obs}}|\Theta)}{\pi (Y_{\text{obs}})}.
\end{align}
Bayes' theorem is an extension of the well-known formula $$P(B|A)=\frac{P(B)P(A|B)}{P(A)}$$
for the conditional probabilities $P(B|A)$ of events $A,B$, to the case of random variables.

The probability $\pi(Y_{\text{obs}})$ in \eqref{posterior-distribution} is unknown, but can be ignored as it is merely a normalizing constant, which in principle could be computed from the fact that
$\pi_{\text{post}}(\Theta)$ is a probability distribution (hence has integral one). The relation in \eqref{posterior-distribution} is usually formulated as
\begin{equation}\label{posterior-distribution_propto}
\pi_{\text{post}}(\Theta) \propto \pi_{\text{prior}}(\Theta) \pi( Y_{\text{obs}}|\Theta).
\end{equation}

\emph{The prior distribution.} In the Bayesian approach, prior knowledge about $\Theta$ is encoded in the prior distribution $\pi_{\text{prior}}(\Theta)$, e.g., from experimental data or known statistical properties, see end of Subsection \ref{Inverse-problem-formulation}\ref{sec: Q to para} for our case.

It is common practice to use a so-called improper or diffusive prior, a constant density on a large interval when there is little knowledge about the prior distribution or when we want to minimise the influence of the prior distribution. In our numerical experiments, we choose (i) the improper prior of the characteristic function of the interval $[0,\infty)$, noting that the parameters $\alpha$ and $\beta$ must be nonnegative, (ii) a Gaussian prior centred at the value of $\Theta$ used for the simulation of the data and truncated at zero. Detailed information and advice on the choice of a prior distribution can be found in \cite[Section 6.3]{HeldBove:2014}.

\emph{The distribution of the error.} The crucial ingredient for constructing the likelihood function is the error $E$, which incorporates uncertainties in the data, the model, and the numerical implementation. The error $E$ is considered to be a mean zero random variable acting as additive noise, and the probability density {$\pi_{\rm err}(E)$}  is assumed to be known. The original problem is modified as follows:
\[
\mathcal{F}(X,\Theta)+E= Y.
\]
For practical implementations, the error term $E = Y_{\rm obs} - \mathcal{F}(X,\Theta)$ is a vector with the same dimension as $Y_{\rm obs}$.
Then, the probability density of the error is taken to be
\begin{equation}\label{error}
   {\pi_{\rm err}(E)} \propto \text{exp}\left( -\frac{1}{2} E^T \Sigma_E^{-1} E \right);
\end{equation} consistent with the assumption that the error $E$ has a Gaussian distribution with zero mean \cite{Kaipio:2005,Tarantola:1987,Yuen:2010}.
The covariance matrix $\Sigma_E$ is usually based on the empirical covariance matrix of the observed data \cite{Rosic:2016} or, simpler, on the identity matrix multiplied by the empirical variance of $Y_{\rm obs}$.

For multiple observation variables, say a bivariate $E = (E_1,E_2)$, the errors are combined to give
\begin{equation}\label{error2D}
   {\pi_{\rm err}(E)} \propto \text{exp}\left( -\frac{1}{2}(E_1^T \Sigma_{E_1}^{-1} E_1 + E_2^T \Sigma_{E_2}^{-1} E_2)\right).
\end{equation}
Let $\sigma_1^2, \sigma_2^2$ be the empirical variances of the components of $Y_{\rm obs}$ and $\Sigma_{E_i}$ (the diagonal matrix with entries $\sigma_i^2$). Then \eqref{error2D} reduces to
\begin{equation}\label{error2Dreduced}
  {\pi_{\rm err}(E)} \propto \text{exp}\left( -\frac{1}{2}(\|E_1\|^2_{\ell^2}/\sigma_1^2 + \|E_2\|^2_{\ell^2}/\sigma_2^2)\right).
\end{equation}
The introduction of the different variances, $\sigma_1^2$ and 
 $\sigma_2^2$ guarantees that differences in scale (magnitude) of the individual components of $Y_{\rm obs}$ are levelled out, so that each component contributes equally to the
error term.

\emph{The likelihood function.}  The likelihood function is obtained by substituting $E = Y_{\rm obs} - \mathcal{F}(X,\Theta)$ into \eqref{error} as shown below:
\begin{equation}\label{likelihood}
\pi(Y_{\text{obs}}|\Theta) = {\pi_{\rm err}}(Y_{\text{obs}} - \mathcal{F}(X,\Theta))
\end{equation}
or equivalently
\begin{align}\label{likelihoodGauss}
  \pi( Y_{\text{obs}}|\Theta) \propto \exp\left( -\frac{1}{2} (Y_{\text{obs}}-\mathcal{F}(X,\Theta) )^T \Sigma_E^{-1} (Y_{\text{obs}}-\mathcal{F}(X,\Theta) ) \right).
\end{align}

The likelihood function is the probability distribution of $Y_{\rm obs}$ when $\Theta$ is given, inherited from the assumed probability distribution of the error between observation $Y_{\rm obs}$ and model output $\mathcal{F}(X,\Theta)$.

\emph{The posterior distribution.} The posterior distribution of $\Theta$ is then given by
\begin{equation}\label{posterior}
    \pi_{\text{post}}(\Theta) \propto \pi_{\text{prior}}(\Theta) {\pi_{\rm err}}(Y_{\text{obs}} - \mathcal{F}(X,\Theta)),
\end{equation}
 obtained by combining \eqref{posterior-distribution_propto} with \eqref{likelihood}. From the Bayesian point of view, it constitutes the solution to the statistical inverse problem \cite[Chapter 3]{Kaipio:2005},
as it provides the probability distribution of the parameter $\Theta$ to be reconstructed, given observation data $Y_{\rm obs}$ and the prior distribution. The factor of proportionality in \eqref{posterior} can be obtained as the reciprocal of the integral
$\int \pi_{\text{prior}}(\theta) {\pi_{\rm err}}(Y_{\text{obs}} - \mathcal{F}(X,\theta)){\rm d}\theta$.
If $\Theta$ is low-dimensional, the factor can be computed by numerical integration, where each step requires an evaluation of the model function $\mathcal{F}(X,\theta)$. However, we shall see in Subsection \ref{Metropolis-algorithm} that the normalization is not required for
sampling from the posterior distribution.

\emph{Bayesian estimators.} Given the posterior distribution, $\pi_{\text{post}}(\Theta)$, one can obtain point estimates $\widehat{\Theta}$ for $\Theta$ \cite[Section 3.1.1]{Kaipio:2005}. We use the mean and median of the posterior distribution as our estimators, $\widehat{\Theta}$, in the subsequent discussion.
The posterior mean is defined as follows: \begin{equation}\label{eq:posterior_mean}
    \widehat{\Theta} = \int \theta\, \pi_{\text{post}}(\theta)\,{\rm d}\theta.
\end{equation}
Given a statistical sample $\theta_1,\ldots,\theta_N$ of the posterior distribution, an approximation to \eqref{eq:posterior_mean} can be obtained as the sample mean $\frac1{N}\sum_{j=1}^N\theta_j$. Then, one can estimate the error as follows.
Suppose the true value to be reconstructed is $\theta^\ast$. Then
\[
   \Big\vert\frac1{N}\sum_{j=1}^N\theta_j - \theta^\ast\Big\vert \leq \Big\vert\frac1{N}\sum_{j=1}^N\theta_j - \int \theta\, \pi_{\text{post}}(\theta)\,{\rm d}\theta\Big\vert
       + \left\vert\int \theta\, \pi_{\text{post}}(\theta)\,{\rm d}\theta - \theta^\ast\right\vert.
\]
The first summand on the right-hand side can be estimated by the usual Monte Carlo error estimate: the average error is of magnitude $\sigma/\sqrt{N}$, where $\sigma^2$ is the variance of $\Theta$. One can do slightly better asymptotically
with quasi-Monte Carlo random generators for which a pointwise error bound ${\mathcal O}(\log N/N)$ can be achieved if $\pi_{\text{post}}(\theta)$ is sufficiently regular \cite{Niederreiter:1992}.
However, there is no control on the second summand, which subsumes a possible systematic error of the approach.

\emph{The observed data.} In our methodological study, $Y_{\rm obs}$ are the numerically computed solutions of \eqref{pde-two-parameters} for fixed values $\alpha_{\ast}$, $\beta_{\ast}$. Given $Y_{\rm obs}$ and $X$ (length scale, temperature), the posterior distributions of $\alpha$, $\beta$ are reconstructed using the Bayesian method outlined above.

\subsection{The Metropolis algorithm -- Markov chain Monte Carlo}
\label{Metropolis-algorithm}

\emph{Markov chains.} We refer to \cite{robert1999monte} for details. A Markov chain is a sequence of random variables $X_0, X_1, X_2,\ldots$ with values in a state space $S$ with the property
-- informally speaking -- that the probability of a transition from one state to another depends only on the current state. More precisely, the conditional probability distributions satisfy
\[
\pi(X_{n+1}|X_0, X_1,\ldots, X_n) = \pi(X_{n+1}|X_n).
\]
We work here solely with a continuous state space, so that all random variables have a probability density. Let $\pi^{(n)}(y)$ be the marginal (i.e., unconditional) probability density of $X_n$.
There is a so-called transition kernel $p(x,y)$, such that at given $x$, $p(x,y)$ is a probability density and at given $y$, $p(x,y)$ is a measurable function of $x$.
The marginal probability densities are computed as
\[
   \pi^{(n)}(y) = \int_S p(x,y)\pi^{(n-1)}(x){\rm d}x
\]
and thus are uniquely determined by the initial density $\pi^{(0)}(x)$.  A probability distribution $\pi$ is called stationary with respect to the Markov chain if
\[
\pi(y) = \int_S p(x,y)\pi(x){\rm d}x.
\]
Under rather general assumptions, a Markov chain has a unique stationary distribution, and the distributions $\pi^{(n)}$ of $X_n$ converge to the stationary distribution $\pi$ in the following sense: Let $\pi^{(0)}(x)$ be any initial density and $f(x)$ be any function integrable with respect to $\pi(x){\rm d}x$. Let $\xi_0,\xi_1,\xi_2,\ldots $ be a realization of the chain. Then
\[
  \lim_{N\to\infty} \frac1{N} \sum_{n=1}^N f(\xi_n) = \int_S f(x)\pi(x){\rm d}x.
\]
This is often referred to as the ergodic theorem \cite[Theorem 6.63]{robert1999monte}. Therefore, the end pieces $\xi_M, \cdots, \xi_N$ for sufficiently large $M$ (after the so-called burn-in phase) are treated as a sample of the limiting distribution $\pi$.
By the ergodic theorem, the sample mean, moments and quantiles of $\xi_M, \cdots, \xi_N$ converge to the expectation value, moments and quantiles of $\pi$
as $N\to\infty$.

\emph{MCMC -- the Metropolis-Hastings algorithm.}
Let $\pi(x)$ be any given probability distribution on a state space $S$. The Metropolis-Hastings algorithm delivers a realization $\xi_0,\xi_1,\xi_2,\ldots $ of a Markov chain which has $\pi(x)$ as its stationary distribution.

Preparation of the algorithm:
\begin{itemize}
\item Choose an initial density $\pi^{(0)}(x)$ or simply an initial point $x = \xi_0$.
\item Choose an arbitrary transition kernel $q(x,y)$ (the so-called proposal distribution) such that $q(x, y) = q(y, x)$
for all $x, y \in S $.
\end{itemize}
Execution of the algorithm:
\begin{itemize}
    \item Sample a value $\xi_0$ from $\pi^{(0)}$.
    \item For $k = 1,\dots, N$, sample a candidate value $\eta$ from the proposal distribution $ q(\xi_{k-1},\cdot)$.
    \item Compute the ratio $r =\frac{\pi(\eta)}{\pi(\xi_{k-1})}.$
    \begin{itemize}
        \item  If $r \geq 1$, the value $\eta$ is accepted; set $\xi_{k}=\eta$.
        \item If $r < 1$, the value $\eta$ is accepted with probability $r$ and rejected
with probability $1-r$.
    \end{itemize}
 \item Draw a random number $\zeta$ from the uniform distribution on $[0, 1]$.
 \begin{itemize}
     \item If $\zeta\leq r,$ set $\xi_k=\eta.$
     \item If $\zeta > r,$ set $\xi_k=\xi_{k-1}.$
 \end{itemize}
\end{itemize}
The algorithm ensures that $\xi_0,\xi_1,\xi_2,\ldots $ is a realization of a Markov chain with $\pi(x)$ as its stationary distribution. In particular, therefore, the end pieces $\xi_M, \cdots, \xi_N$ for sufficiently large $M$ provide a sample of $\pi$.

An important indicator of the quality of the Markov chain is the acceptance rate. This is the proportion of accepted values, that is, $1/N$ times the number of $k$'s for which the candidate value $\eta$ is not rejected. If the acceptance rate is close to one, the chain tends to become the Markov process generated by the proposal distribution alone. If the acceptance rate is small, the chain hardly moves and covers the state space only very slowly. Thus a medium size acceptance rate is desirable. Studies of the efficiency of the Metropolis-Hastings algorithm suggest an optimal acceptance rate of $44\%$ for univariate chains and around $23\%$ for multivariate chains \cite{Gelman:1996}. The convergence of the Markov chain, produced by MCMC, can be proven under very general assumptions on the transition kernel and the stationary distribution, $\pi$
\cite[Section 7.2]{robert1999monte}. We use a transition kernel of the form $q(x,y) = q(x-y)$. In this case, additional requirements on the shape and decay of $\pi$ are needed \cite[Section 7.5]{robert1999monte}, which are usually fulfilled, except for example, in the case of non-identifiable models.

\subsection{Practical aspects of the Metropolis algorithm}
\label{subsec:practicalMCMC}

\begin{figure}[h]
    \centering
    \includegraphics[width = 5.5cm]{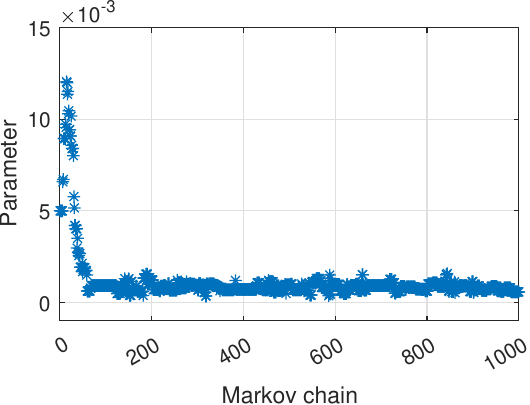} \qquad
    \includegraphics[width = 5.5cm]{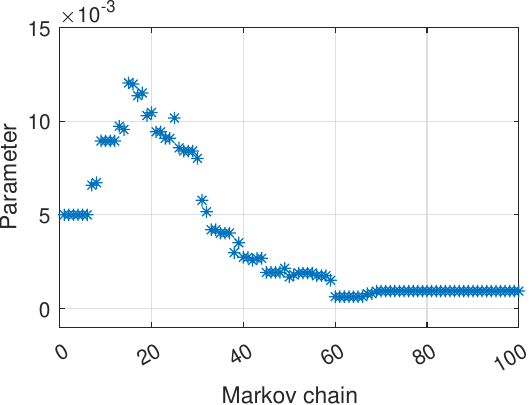}
    \caption{Chain for $\alpha_\ast = 0.0008$ in case of the diagonal solution with uniform prior: first 1000 entries (left), burn-in phase (right).}
    \label{fig:burn-in-phase}
\end{figure}

\emph{Burn-in phase.} The so-called burn-in phase is defined to be the initial segment of the Markov chain, before the elements of the chain hit the main range of the posterior distribution. The length of the burn-in phase depends on the choice of the initial value (and the variance of the proposal distribution); see Figure \ref{fig:burn-in-phase} that corresponds to the $\alpha$-component of the Markov chain for the diagonal solution with $\alpha_\ast = 0.0008$, $\beta_\ast = 1.4$ and a uniform prior distribution.

In the numerical experiments in Sections \ref{numerical-constant-L} and \ref{numerical-two-para}, we simulate a Markov chain of length $10^4$, and we remove the first 200 elements as burn-in phase for the statistical analysis.

\emph{Monitoring convergence.}
The quality of the simulated Markov chain can be tested by means of its mixing properties, acceptance rate, and tests for stationarity, 
Additional methods of analysing the quality of Markov chains produced by MCMC can be found in
\cite[Chapter 12]{robert1999monte} and in \cite{Brooks:1998}. For further  statistical explorations of posterior distributions see, e.g., \cite{Barndorff-Nielsen:1994}
and \cite{Efron:2015}.

To assess the quality of the Markov chain, one usually checks if the chain explores the whole state space (the range of the parameter under consideration) reasonably fast. One important criterion is the acceptance rate, which optimally should be around 44\% for univariate chains around $23\%$ for multivariate chains, see Subsection \ref{Metropolis-algorithm}. The acceptance rates obtained in our numerical experiments are recorded in Tables 
\ref{tab:statsD1}, \ref{tab:statsR4},
\ref{tab:statsD1R4largealpha} and \ref{tab:statsD1R4smallalpha} and are reasonably close to the optimal ones. A second criterion is a good mixing property so that the candidates move between different parts of the state space quickly. This can be checked visually and is evident for all chains in the figures in Sections~\ref{numerical-constant-L} and \ref{numerical-two-para}.

The final question concerns the convergence of the chain to its stationary distribution. This means that subsamples $\xi_K, \ldots, \xi_{K+M}$ and $\xi_L, \ldots, \xi_{L+M}$, with sufficiently large $K$ and $L>K+M$, should have the same distribution. This can be tested by nonparametric statistical tests such as the Kolmogorov-Smirnov test. However, these tests apply only to iid (independent and identically distributed) subsamples.  As suggested in \cite[Section 12.2.2]{robert1999monte}, independent sampling can be achieved by selecting batches of the form
$\xi_K, \xi_{K + G}, \xi_{K+2G},\ldots, \xi_{K+lG}$ with $G >\!\!> 1$ (and $lG\leq M$). We perform these tests repeatedly for successive periods of length 1000 with batch step $G = 10$. The obtained values of the Kolmogorov-Smirnov test statistic (for two samples of equal size) are well below the critical value \cite[Table ST8]{Rohatgi:2015}
at significance level 5\%, in most cases.

\emph{Confidence intervals}
The approximate confidence intervals of the mean in Section \ref{numerical-constant-L} (Figure \ref{fig:stat_vs_length}) can be computed from a Central Limit Theorem for Markov chains \cite[Theorem 27.4]{Billingsley:1995}: Let $X_1, X_2, \ldots$ be a stationary Markov chain satisfying certain mixing conditions; let $S_N = X_1 + \cdots + X_N$. Then the variance ${\rm V}(S_N/N)$ of $S_N/N$ converges to
\begin{equation}\label{eq:gamma2}
   \gamma^2 = {\rm V}(X_1) + 2\sum_{k=1}^\infty{\rm COV}(X_1,X_{1+k})
\end{equation}
as $N\to\infty$. If $\gamma > 0$, the random variable $S_N/\gamma\sqrt{N}$ converges (in distribution) to the standard normal distribution. The required mixing conditions
cannot easily be verified in practice, but can be enforced in the MCMC setting (see the discussion in \cite[Sections 6.7.2 and 6.9]{robert1999monte}).

Thus, taking a segment $\xi_{M+1}, \ldots, \xi_{M+N}$ of length $N$ of the chain after the burn-in period with sample mean $\mu_{\alpha,N}= \frac{1}{N}\sum_{i=1}^N\xi_{M+i}$, one may assume that
$\mu_{\alpha,N}$ is approximately normally distributed around $\mu^\ast$ with variance $\gamma^2/N$; here $\mu^\ast$ is the true mean of the posterior distribution.
Due to the assumed stationarity of the chain, the variance and the covariances in \eqref{eq:gamma2} can be estimated by the empirical autocovariance function of the segment $\xi_{M+1}, \ldots, \xi_{M+N}$, and the infinite sum can be truncated at large $k = K$ when the covariances become negligible (empirically around $K = 10$ to $20$).

Based on these considerations, approximate confidence intervals for the posterior mean can be computed as $\mu_{\alpha,N} \pm q\gamma/\sqrt{N}$, where $q$ is the quantile of the standard normal distribution corresponding to the desired confidence level ($q=1.96$ for 95\% confidence), as has been done in Figure \ref{fig:stat_vs_length}.

\section{Computational framework} \label{sec:computational_framework}

In this section, we briefly describe our finite element framework for solving the system in \eqref{pde-two-parameters}, subject to prescribed, continuous and piecewise smooth Dirichlet boundary data, $\Qvec_b$. We take $\Omega$ to be a convex polygon and let $\mathcal{T}$ be a shape regular triangulation of $\Omega$ into triangles ($T$). The maximal mesh size in a triangulation is denoted by $h:=\max_{T\in \mathcal{T}} \text{diam}(T)$, where $\text{diam}(T)$ is the diameter of each simplex $T\in \mathcal{T}.$ Define the finite element space $X_h:= \{\varphi \in C^0(\Bar{\Omega})| \, v|_T \in P_1(T) \text{ for all } T\in \mathcal{T} \} \subset H^1(\Omega)$ and  $X_h^0:=X_h\cap  H^1_0(\Omega),$ where $P_1(T)$ is the space of all polynomials of degree at most $1$ on $T.$ The Galerkin discretization of the nonlinear problem \eqref{pde-two-parameters} is to find $\Qvec_h:= (Q_{11,h}, Q_{12,h}) \in X_h \times X_h$ satisfying the Dirichlet boundary condition $\Qvec_h|_{\partial \Omega}:=\mathcal{I}_h \Qvec_b$, obtained by conforming interpolation, such that
\begin{align}\label{equn:discrete-PDE}
\begin{split}
    \alpha  \int_\Omega \nabla Q_{11,h} \cdot \nabla \varphi_{1,h} \dx+ \int_\Omega (Q_{11,h}^2+Q_{12,h}^2 -\beta)Q_{11,h} \varphi_{1,h}\dx =0,
  \\
  \alpha  \int_\Omega \nabla Q_{12,h} \cdot \nabla \varphi_{2,h} \dx+ \int_\Omega (Q_{11,h}^2+Q_{12,h}^2 -\beta)Q_{12,h} \varphi_{2,h}\dx =0 ,
\end{split}
\end{align}
for all $\Phi_h:=(\varphi_{1,h},\varphi_{2,h}) \in X_h^0 \times X_h^0 $.

We use Newton's method to approximate solutions of the discrete nonlinear system \eqref{equn:discrete-PDE}.
The initial guess plays a crucial role in the convergence of the Newton iterates and the selected solution. We abstractly write the  nonlinear system  \eqref{equn:discrete-PDE} as an operator equation $\mathcal{N}_h(\Qvec_h )=0$. Starting from
an initial guess $\Qvec^0_h$, Newton iterates $\Qvec^{k+1}_h=\Qvec^k_h+\delta\Qvec$ are computed by solving $\langle D \mathcal{N}_h(\Qvec^k_h)\delta\Qvec,  \Phi_h\rangle = -\mathcal{N}_h (\Qvec^k_h, \Phi_h)$. Here, $D \mathcal{N}_h(\Qvec^k_h)$ denotes the Fr\'{e}chet derivative of $\mathcal{N}_h$ at $\Qvec^k_h$. We refer to  \cite{DGFEM,AbstractAMRMNN2023} for more algorithmic details.

\subsection*{Benchmark example }

We review the benchmark example of NLCs confined to a square domain with tangent boundary conditions, as motivated by the planar bistable nematic device reported in \cite{Tsakonas} and the work in \cite{MultistabilityApalachong}. 
Consider the rescaled square domain $\Omega:=[0,1] \times [0,1] $ with tangent boundary conditions on the square edges, in the reduced LdG framework outlined in Subsection \ref{Inverse-problem-formulation}\ref{sec: Q to para}. The tangent boundary conditions mean that on the square edges, the nematic director has to be tangent to the edges, i.e., $\nvec = (\pm 1, 0)$ on the horizontal edges $y=0,1$, and $\nvec = (0, \pm 1)$ on the vertical edges, $x=0,1$. Fixing $s=1$ on the square edges (this is an arbitrary choice that essentially means that the molecules are well ordered on the edges), this translates to the boundary conditions displayed in Table \ref{table:tangential-BC}.  The discontinuities at the vertices are circumvented by defining the Dirichlet boundary condition, $\Qvec_b$, using the trapezoidal shape function \cite{MultistabilityApalachong}; see below:
	\begin{equation} \label{BC}
	\Qvec_b(x,y)=
	\begin{cases}
	(\textit{T}_{d}(x),0) & \text{on}\,\,\,\, y=0 \,\,\,\,\text{and} \,\,\,\, y=1, \\
	(- \textit{T}_d(y),0)  & \text{on}\,\,\,\, x=0 \,\,\,\,  \text{and} \,\,\,\,  x=1,
	\end{cases}
	\end{equation}
	with parameter $d=0.06$  and the trapezoidal shape function $\textit{T}_d:[0,1]\rightarrow {\mathbb{R}}$ is defined by
	\begin{equation*}
	\textit{T}_d(t)=
	\begin{cases}
	t/d, & 0 \leq t \leq d,  \\
	1, &  d \leq  t \leq 1- d, \\
	(1-t)/d, & 1- d \leq t \leq 1.
	\end{cases}
	\end{equation*}	
\begin{tiny}
\begin{table}
\centering
\begin{tabular}{|c | c | c | c | c |}
\hline
Solution & $x=0$ & $x=1$ & $y=0$ & $y=1$\\
\hline
$Q_{11}$ & $-1$ & $-1$ & 1 & 1 \\
$Q_{12}$ & $0$ & $0$ & $0$ & $0$ \\
\hline
\end{tabular}
\caption{Tangent boundary conditions  for $\Qvec.$}
\label{table:tangential-BC}
\end{table}
\end{tiny}

\begin{figure}[h]
    \centering
    \subfloat[Mesh]{\includegraphics[width=0.28\linewidth]{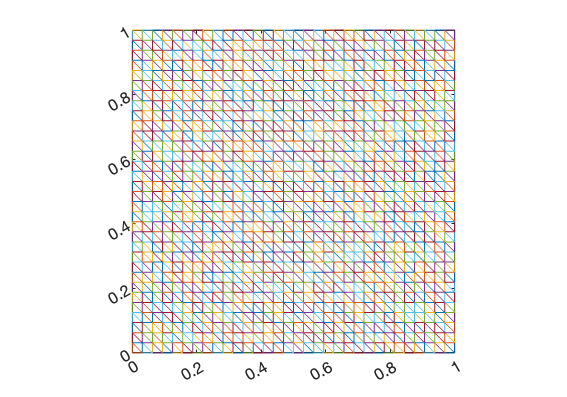}}
    \subfloat[Diagonal solution]{\includegraphics[width=0.28\linewidth]{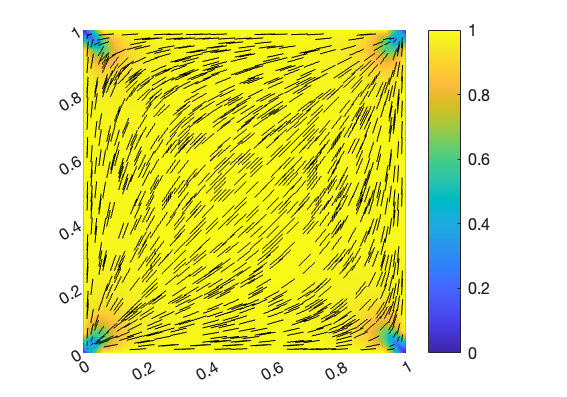}}
     \subfloat[Rotated solution]{\includegraphics[width=0.28\linewidth]{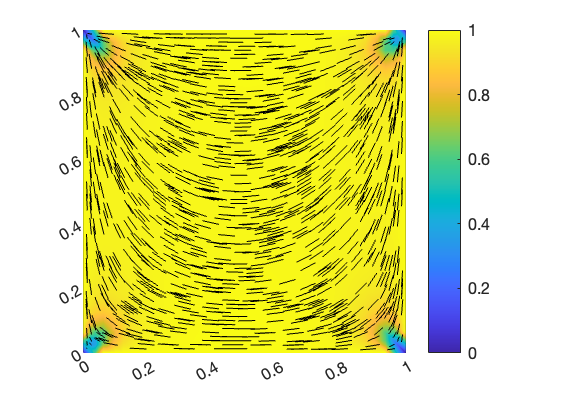}}
    \caption{Mesh  with $h=0.0442$;  the diagonal and rotated solutions computed on this mesh. 
    }
    \label{fig:mesh-diagonal-soln}
\end{figure}

It is known both experimentally and numerically that the  system \eqref{equn:discrete-PDE}, subject to the Dirichlet boundary condition $\Qvec_b$ defined above, admits two classes of solutions for small enough $\alpha$ (for $\lambda \geq \lambda_c$ where $\lambda_c$ is estimated to be in the nanometre range from the results in \cite{Karlj_Majumdar2014}):
\begin{enumerate}
\item diagonal solutions: the planar director, modelled by $\nvec$ in (\ref{eq:Q}), roughly aligns along one of the square diagonals and there are two classes of diagonal solutions: D1 and D2, one for each square diagonal;
\item rotated solutions: $\nvec$ rotates by $\pi$ radians between a pair of opposite square edges, and there are 4 classes of rotated solutions labelled by R1, R2, R3 and R4 respectively, related to each other by $\pi/2$ radians.
\end{enumerate}
In other words, there are at least six competing solutions (two diagonal and four rotated) of the system (\ref{pde-two-parameters}) for $\lambda$ large enough, or $\alpha$ small enough.
In Figure \ref{fig:mesh-diagonal-soln}, we plot the computational mesh for the rescaled unit square domain, and the numerically computed diagonal and rotated solutions. The diagonal and rotated solutions are distinguished by the locations of the splay vertices; a splay vertex being a vertex such that $\nvec$ splays around it. 
Each diagonal solution has a pair of diagonally opposite splay vertices and each rotated solution has a pair of adjacent splay vertices, connected by an edge of the square. The black arrows correspond to $\mathbf{n}=(\cos \theta,\sin \theta )$ in (\ref{eq:Q}) with $\theta:=\frac{1}{2} \arctan \frac{Q_{12}}{Q_{11}}$  whilst the colour bar refers to the scalar order parameter $s$ in (\ref{eq:Q}).

\section{Numerical experiments: one parameter identification}\label{numerical-constant-L}

Next, we demonstrate how to reconstruct the parameters $\alpha,\beta$ in \eqref{pde-two-parameters} from given $\Qvec$, using the Bayesian methods and the Metropolis-Hastings algorithm outlined in Section~\ref{sec:Bayes}. We focus on the benchmark square domain example in Section~\ref{sec:computational_framework} and perform numerical experiments for various values of $\alpha \in \left[0.0008, 0.004 \right]$ and $\beta \in \left[0.6, 1.4\right]$; cf.~end of Subsection \ref{Inverse-problem-formulation}\ref{sec: Q to para}.%

 In this section, we fix $\beta=1$ and focus on the reconstruction of $\alpha$ in \eqref{pde-two-parameters}, and in the next section, we reconstruct both $\alpha$ and $\beta$ for given (measured) $\Qvec$.

The strategy is as follows. We first generate artificial data $\Qvec_{\text{obs}}$ by numerically solving the nonlinear system \eqref{equn:discrete-PDE} on a square domain with the boundary conditions in Table~\ref{table:tangential-BC}, with reference values of $\alpha = \alpha_{\ast}$ and $\beta=1$.
 The full solution of the statistical inverse problem is then given by the posterior distribution of $\alpha$ as outlined in Section\;\ref{sec:Bayes}. From there, the posterior mean or the posterior median serves as point estimate $\widehat{\alpha}$. More precisely, substituting $\Theta = \alpha$ and $Y_{\text{obs}} = \Qvec_{\text{obs}}$ in formula \eqref{posterior} yields
\begin{equation*}
    \pi_{\text{post}}(\alpha) \propto \pi_{\text{prior}}(\alpha) {\pi_{\rm err}}(\Qvec_{\text{obs}} - \mathcal{F}(X,\alpha)).
\end{equation*}
Here, the variable $X$ contains information about the square edge length and ambient temperature and $\mathcal{F}(X,\alpha)$ is the solution of the discretized forward problem \eqref{equn:discrete-PDE} corresponding to a fixed value of $\alpha$
(and $\beta = 1$). Recalling that $\Qvec=(Q_{11},Q_{12})$ and $\mathcal{F}=(\mathcal{F}_{11},\mathcal{F}_{12})$ have two components, and assuming a Gaussian error distribution \eqref{error2Dreduced}, \eqref{likelihoodGauss} yields
\begin{align}\label{likelihood-LC-example}
  \pi_{\text{post}}( \alpha) \propto \pi_{\text{prior}}(\alpha)\exp\left( -\frac{1}{2} \left(\norm{\bar{Q}_{11}-\mathcal{F}_{11}(X,\alpha)}_{\ell^2}^2
   /\sigma_{11}^2 + \norm{\bar{Q}_{12}-\mathcal{F}_{12}(X,\alpha)}_{\ell^2}^2/\sigma_{12}^2 \right) \right),
\end{align}
where $\Qvec_{\text{obs}}=(\bar{Q}_{11},\bar{Q}_{12})$ are the generated reference data and $\sigma_{11}^2,$ $\sigma_{12}^2$ are the corresponding variances (computed empirically from the values of $\bar{Q}_{11},\bar{Q}_{12}$ at the grid points of the triangulation  $\mathcal{T}$).

We employ two choices of the prior distribution for $\alpha$: (i) an improper uniform prior, and (ii) a Gaussian prior distribution.
The improper uniform prior (UP) distribution is chosen to be $\pi_{\text{prior}}(\alpha)  = \chi(\alpha)$, where $\chi$ is the characteristic function of $(0,\infty)$ which encodes that $\alpha$ is positive. The Gaussian prior (GP) is chosen to be $\pi_{\text{prior}}(\alpha)  \propto \chi(\alpha) \text{exp}\big( -(\alpha-\alpha_{\ast})^2/2\sigma_0^2 \big)$ where $\alpha_{\ast}$ is the reference parameter used for generating the artificial data set $\Qvec_{\text{obs}}$. We select the standard deviation to be $\sigma_0=0.0005$ for all our test cases. This choice guarantees a reasonable coefficient of variation, $\sigma_0/\alpha_\ast$, of 12.5\% for $\alpha_{\ast} = 0.004$. 
To implement the Metropolis-Hastings algorithm, we choose a proposal distribution of the form  $q(x,y) = q(x-y)$, where $q(\bullet)$ is the density of a mean-zero Gaussian distribution with
standard deviation $0.001$. This guarantees good mixing properties of the generated Markov chain  (see Section\;\ref{sec:Bayes}\ref{subsec:practicalMCMC}, a good mixing property is usually guaranteed when the standard deviation is of the same or somewhat smaller magnitude than the parameter to be reconstructed.). The initial value of the Markov chain is always chosen to be $0.005$, which results in a rather short burn-in phase for all cases. 

Next, we use the diagonal and rotated solutions of \eqref{equn:discrete-PDE} (for the benchmark example in Section~\ref{sec:computational_framework}) with $\alpha_\ast=0.004$ and $\beta_\ast=1$ as the reference solution and use the inverse problem approach in Section~\ref{sec:Bayes} to reconstruct the value of $\alpha$. For each choice of the prior distribution of $\alpha$ (UP or GP), we plot the Markov chain, the resulting histogram of the posterior distribution of $\alpha$, tabulate the mean and median as estimates for the reconstructed value of $\alpha$ along with the acceptance rate of the Markov chain which serves as an indicator of the quality of the chain. 

\begin{figure}[h]
    \centering
    \subfloat[Reference solution ]{\includegraphics[width=0.28\linewidth]{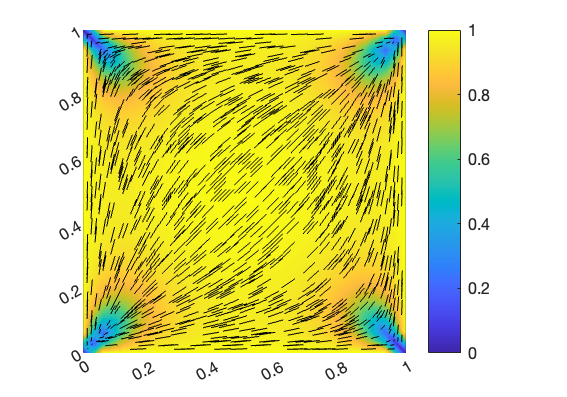}}
    \subfloat[Markov chain for $\alpha$ (UP)]{\includegraphics[width=4 cm, height = 3.2 cm]{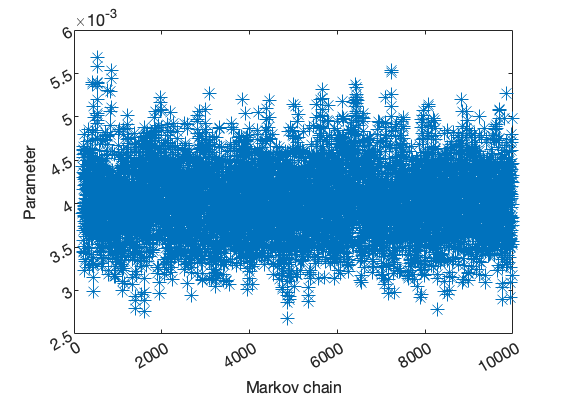}}
     \subfloat[Histogram of $\alpha$ (UP)]{\includegraphics[width=4 cm, height = 3.2 cm]{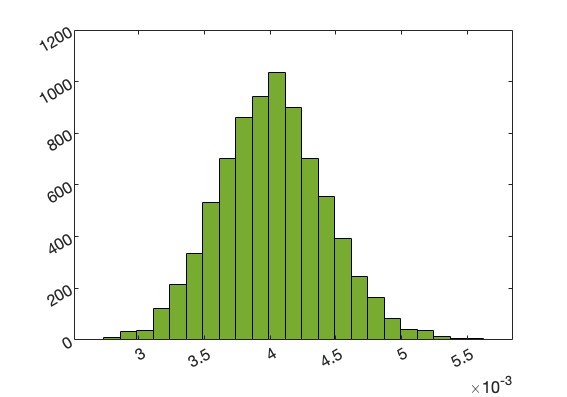}}\\
     \subfloat[Markov chain for $\alpha$ (GP)]{\includegraphics[width=4 cm, height = 3.2 cm]{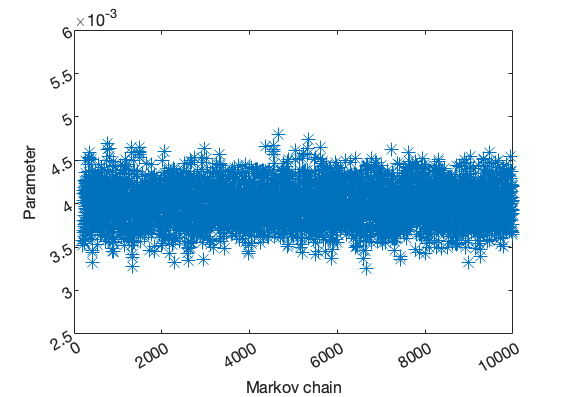}}
     \subfloat[Histogram of $\alpha$ (GP)]{\includegraphics[width=4 cm, height = 3.2 cm]{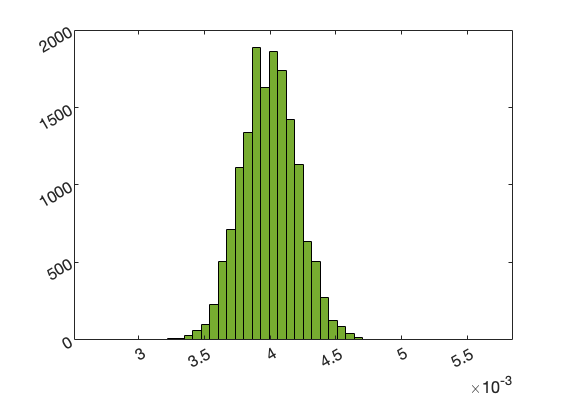}}
\caption{Diagonal solution, reference value $\alpha_{\ast}=0.004$: (a) plot of reference solution; (b) Markov chain and (c) histogram
  of the posterior distribution of parameter $\alpha$ for uniform prior (UP) distribution; (d) Markov chain and (e) histogram
  of the posterior distribution of parameter $\alpha$ for Gaussian prior (GP) distribution. 
}
   
    \label{fig:L-0.008-diag}
\end{figure}

The diagonal solution of \eqref{equn:discrete-PDE} is plotted in Figure~\ref{fig:L-0.008-diag}(a), with $\alpha_\ast=0.004$ and $\beta_\ast=1$, using the mesh plotted in Figure~\ref{fig:mesh-diagonal-soln}(a). 
This solution generates the data $\Qvec_{\text{obs}}$.   Figure  \ref{fig:L-0.008-diag}(b) displays the end piece of length $9800$ of the Markov chain  (with burn-in phase of first $200$ entries removed) obtained by the MCMC algorithm, with the posterior density distribution given by \eqref{likelihood-LC-example}. In Figure  \ref{fig:L-0.008-diag}(c), we plot the histogram associated with this chain, for the uniform prior distribution. In Figures  \ref{fig:L-0.008-diag}(d) and (e), we plot the Markov chain and the histogram associated with the Gaussian prior distribution for $\alpha$ respectively.

The corresponding values of the sample
mean $(\mu_\alpha)$ and the sample median $(m_\alpha)$ are summarized in Table \ref{tab:statsD1}. 
\begin{table}[h]
\centering
\begin{tabular}{|c|c|c|c|c|c|}
\hline
Statistics for Figure \ref{fig:L-0.008-diag} & mean $\mu_\alpha$ & median $m_\alpha$ & standard deviation 
 & acceptance rate
\\ \hline
Uniform prior  & 0.0040195 &   0.0040126 &   0.0004108 
 &   65\%\\
Gaussian prior &  0.0039962  &  0.0039969 &  0.0002121 
 &  45\% \\
\hline
\end{tabular}
\caption{Statistics for diagonal solution for $\alpha_{\ast} = 0.004$.}
\label{tab:statsD1}
\end{table}

Firstly, the posterior mean overestimates the true value in the case of a uniform prior; it is biased. (This is generally the case as can be seen from explicit calculations comparing Bayesian with frequentist estimates as, e.g.,~in \cite[Example 6.8]{HeldBove:2014}.)
Secondly, the posterior median is less biased, as has also been empirically confirmed in the literature \cite{Pick:2023}.

\begin{figure}[h]
    \centering
    \includegraphics[width = 6.4cm]{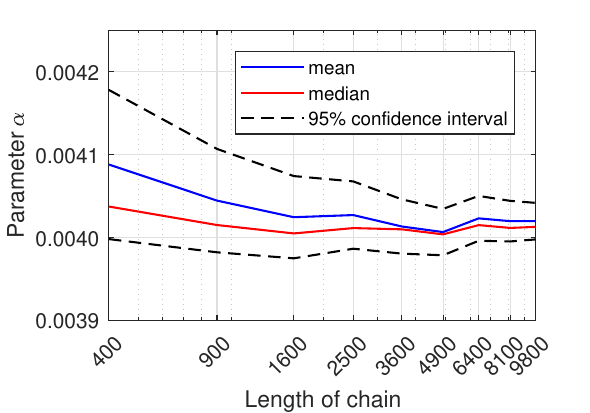}
    \caption{Chain for $\alpha_{\ast} = 0.004$, diagonal, uniform prior: mean, median and 95\% confidence interval for mean depending on length of chain segment.}
    \label{fig:stat_vs_length}
\end{figure}

Asymptotic confidence intervals for the mean can be computed as explained in Subsection \ref{sec:Bayes}\ref{subsec:practicalMCMC}. We also study the trends in the mean, the median and the confidence intervals by taking segments of increasing length $N$ of the Markov chain starting at $M = 201$. The results are summarized 
in Figure\;\ref{fig:stat_vs_length}. One can observe a decreasing bias of the mean and the median, an improved accuracy indicated by a smaller confidence interval and the fact that the reference value $\alpha_{\ast} = 0.004$ always lies in the 95\%-confidence interval.

\begin{figure}[h]
    \centering
    \subfloat[Reference solution ]{\includegraphics[width=0.28\linewidth]{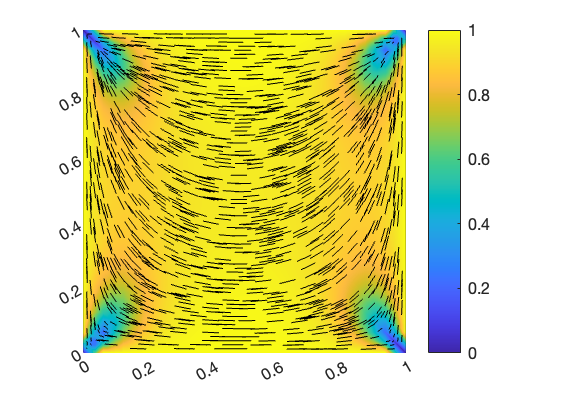}}
    \subfloat[Markov chain for $\alpha$ (UP)]{\includegraphics[width=4 cm, height = 3.2 cm]{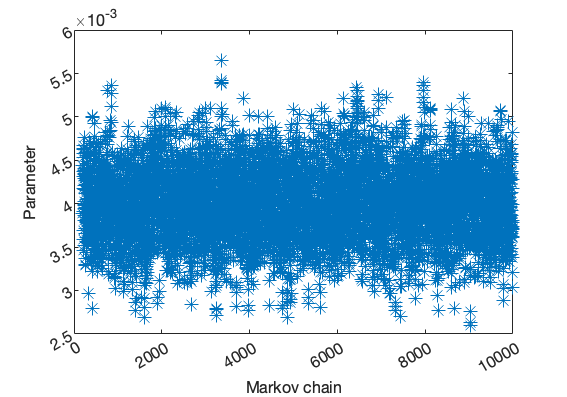}}
     \subfloat[Histogram of $\alpha$ (UP)]{\includegraphics[width=4 cm, height = 3.2 cm]{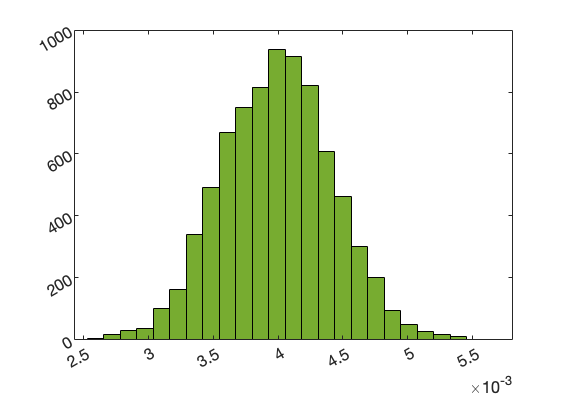}}\\
     \subfloat[Markov chain for $\alpha$ (GP)]{\includegraphics[width=4 cm, height = 3.2 cm]{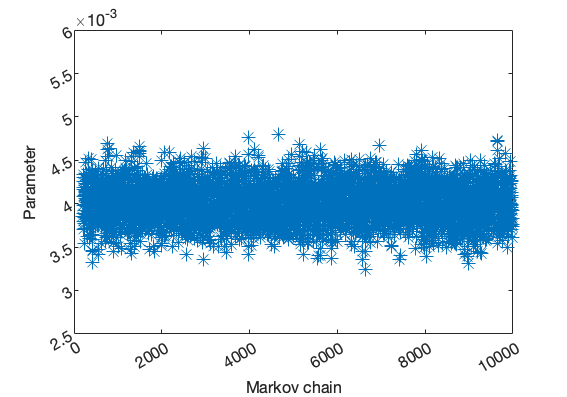}}
     \subfloat[Histogram of $\alpha$ (GP)]{\includegraphics[width=4 cm, height = 3.2 cm]{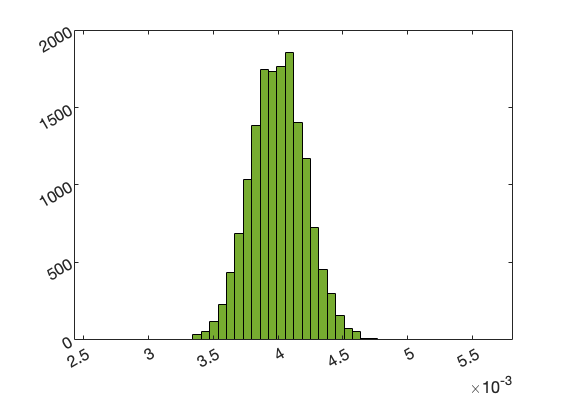}}
\caption{
Rotated solution, reference value $\alpha_{\ast}=0.004$: (a) plot of reference solution; (b) Markov chain and (c) histogram
  of the posterior distribution of parameter $\alpha$ for uniform prior (UP) distribution; (d) Markov chain and (e) histogram
  of the posterior distribution of parameter $\alpha$ for Gaussian prior (GP) distribution. 
}
    \label{fig:L-0.008-rota}
\end{figure}

Next, we repeat the same numerical experiments with the rotated solution of \eqref{equn:discrete-PDE}, with $\alpha_{\ast} = 0.004$ and $\beta_{\ast}=1$ as the reference $\Qvec_{\text{obs}}$ (see Figure \ref{fig:mesh-diagonal-soln}(c)). It is important to test the inverse problem approach in Section~\ref{sec:Bayes} with different solutions, not only to assess the robustness of the method but also to understand whether some solutions outperform others for our purposes and if so, why. 
The description of Figure\;\ref{fig:L-0.008-rota}
is analogous to the explanation of Figure\;\ref{fig:L-0.008-diag}.
The corresponding mean and median values of the posterior distribution of $\alpha$ are summarized in Table \ref{tab:statsR4}.
\begin{table}[h]
\centering
\begin{tabular}{|c|c|c|c|c|c|c|}
\hline
Statistics for Figure \ref{fig:L-0.008-rota} & mean $\mu_\alpha$ & median $m_\alpha$ & standard deviation & acceptance rate\\ \hline
Uniform prior  & 0.0039951 &   0.0039996 &   0.0004258 &   66\%\\
Gaussian prior &  0.0040011  &  0.0039999 &  0.0002129 &  44\% \\\hline
\end{tabular}
\caption{Statistics for rotated solution for $\alpha_{\ast} = 0.004$.}
\label{tab:statsR4}
\end{table}
There are no discernible differences between the rotated and diagonal solutions, or the uniform and Gaussian prior distributions for $\alpha$, but these numerical experiments are not exhaustive and are limited in many ways. 

\section{Numerical experiments: Identification of $\alpha$ and $\beta$ from the diagonal and rotated solutions}
\label{numerical-two-para}

This section demonstrates how the two parameters $(\alpha,\beta)$ in \eqref{pde-two-parameters} can be reconstructed simultaneously. We follow the same strategy as described at the beginning of Section \ref{numerical-constant-L}, using solutions to \eqref{pde-two-parameters} with chosen reference parameters,
$\alpha = \alpha_{\ast}$ and $\beta = \beta_{\ast}$, on a square domain with tangent boundary conditions in Table~\ref{table:tangential-BC}, as artificial data $\Qvec_{\text{obs}}$. The likelihood function is as in \eqref{likelihood-LC-example} and the model outputs
$\mathcal{F}_{11}=\mathcal{F}_{11}(X,\alpha,\beta)$ and $\mathcal{F}_{12}=\mathcal{F}_{12}(X,\alpha,\beta)$ depend on the two parameters $(\alpha,\beta)$.
As prior distributions for $\alpha$ and $\beta$, we take (i) an improper uniform prior $\pi_{\text{prior}} = \chi(\alpha)\chi(\beta)$, where $\chi$ is the characteristic function of the interval $(0,\infty)$ and (ii) a bivariate Gaussian prior as described in the subsequent examples. We consider reference values $\alpha_{\ast} = 0.004$, $\beta_{\ast} = 0.6$ and $\alpha_{\ast} = 0.0008$, $\beta_{\ast} = 1.4$, thus covering a relatively large and a relatively small value of $\alpha_{\ast}$ respectively.
The computations are done with the diagonal and rotated solutions as before, with both prior distributions for $\alpha$ and $\beta$ respectively. The Gaussian prior is always chosen to be centred at $\alpha_\ast$ with standard deviations $\sigma_\alpha=0.0005$, $\sigma_\beta=0.1$, and correlation coefficient $\rho_{\alpha\beta}=\frac{\sigma_{\alpha \beta}}{\sigma_\alpha \sigma_\beta}=0.5$, where $\sigma_{\alpha \beta}$ denotes the covariance of the bivariate Gaussian.

In all cases, the Markov chain is based on a bivariate mean zero Gaussian distribution as proposal distribution $q(\bullet)$ with correlation coefficient $0.8$. For $\alpha_{\ast} = 0.004$, $\beta_{\ast} = 0.6$, the standard deviations of the proposal distribution are fixed to be $(0.001,0.1)$; for $\alpha_{\ast} = 0.0008$, $\beta_{\ast} = 1.4$, the standard deviations are fixed to be $(0.005,0.1)$. The initial points of the Markov chains are $(0.01,0.5)$ for $(\alpha_{\ast}, \beta_{\ast}) = (0.004, 0.6)$,  and $(0.005,0.8)$ for $\left(\alpha_{\ast}, \beta_{\ast}\right) = (0.0008,1.4)$ respectively.
As in Section 5, a Markov chain of length 10000 is generated and the end pieces of length 9800 are used for the analysis. 

Consider a diagonal and a rotated solution of \eqref{pde-two-parameters} with $(\alpha_{\ast}, \beta_{\ast}) = (0.004, 0.6)$ and $\left(\alpha_{\ast}, \beta_{\ast}\right) = (0.0008,1.4)$ respectively, on a square domain with tangent boundary conditions as outlined in Table~\ref{table:tangential-BC}, both of which are treated as $\Qvec_{\text{obs}}$ in our numerical experiments. In what follows, we tabulate the mean and median estimators for $\alpha$ and $\beta$, for both pairs of reference solutions and reference values of $(\alpha_{\ast}, \beta_{\ast})$, and for both choices of the prior distributions for $\alpha$ and $\beta$. The figures are  organised as follows: (a) illustrates the reference solution. A plot of the Markov chains (with burn-in phase deleted) for $\alpha$ and $\beta$
are displayed in (b) and (d), respectively. Plots (c) and (e) focus on the histograms of the (marginal) posterior distributions of $\alpha$ and $\beta$ respectively and (f) is a contour plot of the bivariate histogram of the joint posterior distribution of $(\alpha,\beta)$.

The results are summarised in Tables~\ref{tab:statsD1R4largealpha} and \ref{tab:statsD1R4smallalpha} respectively. D1 labels a diagonal solution and R4 labels a rotated solution respectively. As before, there are no discernible differences between the diagonal and rotated solutions, and between the UP and GP, for $\alpha_{\ast}=0.004$ and $\beta_{\ast}=0.6$. The acceptance rate for the rotated solution is typically higher than for the diagonal solution. 
The estimates of the parameter $\alpha$ in Table \ref{tab:statsD1R4largealpha} are slightly less accurate than in the one-parameter case (Tables \ref{tab:statsD1} and \ref{tab:statsR4}), but still with a relative error between 1\% and 5\%. The relative error for $\beta$  is below 1\%.
However, the rotated solution provides poorer estimates for $\alpha$, when $\alpha_{\ast} = 0.0008$ and $\beta_{\ast} = 1.4$, compared to the diagonal solution, showing a bias of the mean of about $+16\%$. 
We recall that a smaller value of $\alpha_{\ast}$ corresponds to a larger square domain and diagonal solutions are global energy minimisers on large square domains \cite{MultistabilityApalachong}. The rotated solutions are locally stable on large square domains and the energy difference between the diagonal and rotated solutions increases as the square edge length increases. This could provide some heuristic insight into why diagonal solutions perform better than rotated solutions for parameter inference or inverse problems, as the square edge length increases or equivalently as $\alpha_{\ast}$ decreases. 
The estimates of the larger parameter $\beta$ are more accurate with a relative error of less than 0.5\% for the diagonal solution and less than 1\% for the rotated solution.

\begin{table}[h]
\centering
\begin{tabular}{|c|c|c|c|c|c|c|}
\hline
Type of chain & mean $\mu_\alpha$  & median $m_\alpha$ &  mean $\mu_\beta$ & median $\mu_\beta$ & correlation $\rho_{\alpha\beta}$ & AR\\ \hline
D1, UP & 0.0041841	& 0.0041460	& 0.6060851 &	0.6054822 &	0.8654732 &	17\%\\
D1, GP & 0.0040616	& 0.0040508	& 0.6017940	& 0.6020828	& 0.8332836	& 14\%\\
R4, UP & 0.0041191 &	0.0041147 &	0.6085217 &	0.6077666	 &0.9258516 &	22\%\\
R4, GP & 0.0040311 &	0.0040217	 &0.6023924 &	0.6016717 &	0.9064120 &	19\%\\\hline
\end{tabular}
\caption{Sample statistics for posterior distributions, reference values $\alpha_{\ast} = 0.004, \beta_{\ast} = 0.6$.}
\label{tab:statsD1R4largealpha}
\end{table}

\vspace*{-1.2cm}

\begin{table}[h]
\centering
\begin{tabular}{|c|c|c|c|c|c|c|}
\hline
Type of chain & mean $\mu_\alpha$  & median $m_\alpha$ &  mean $\mu_\beta$ & median $\mu_\beta$ & correlation $\rho_{\alpha\beta}$ & AR \\ \hline
D1, UP & 0.0008587 &	0.0008478	& 1.4044713	& 1.4055363 &	0.4372786 &	16\%\\
D1, GP & 0.0008596 &	0.0008417	& 1.4037591 &	1.4046242 &	0.4318876	& 15\%\\
R4, UP & 0.0009278 &	0.0008725 &	1.4134830 &	1.4139106 &	0.4805958	& 31\%\\
R4, GP & 0.0009292 &	0.0008886	& 1.4113185 &	1.4106873 &	0.4865024 &	30\%\\\hline
\end{tabular}
\caption{Sample statistics for posterior distributions, reference values $\alpha_{\ast} = 0.0008, \beta_{\ast} = 1.4$.}
\label{tab:statsD1R4smallalpha}
\end{table}

\vspace*{-0.7cm}

\begin{figure}[H]
    \centering
      \subfloat[Reference solution]{ \includegraphics[width=0.28\linewidth]{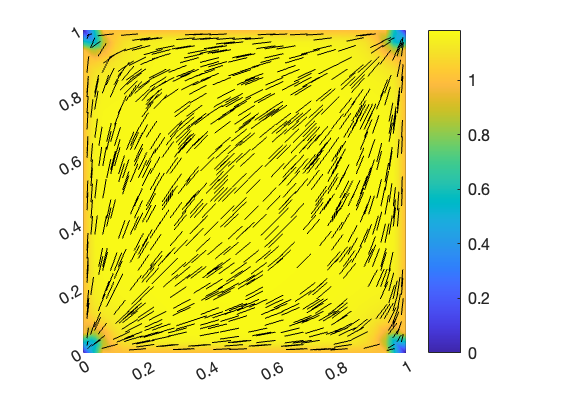}}
    \subfloat[Markov chain for $\alpha$]{\includegraphics[width=4 cm, height = 3.2 cm]{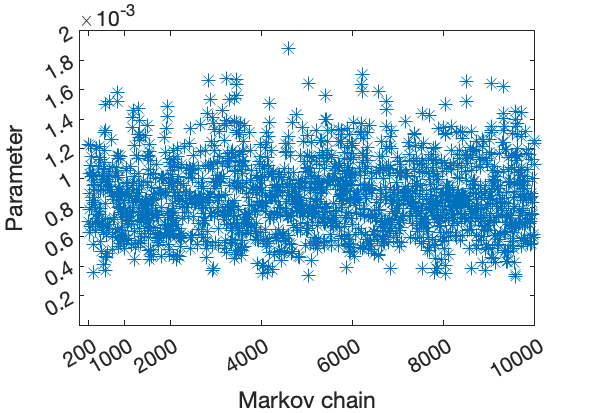}}
     \subfloat[Histogram of $\alpha$]{\includegraphics[width=4 cm, height = 3.2 cm]{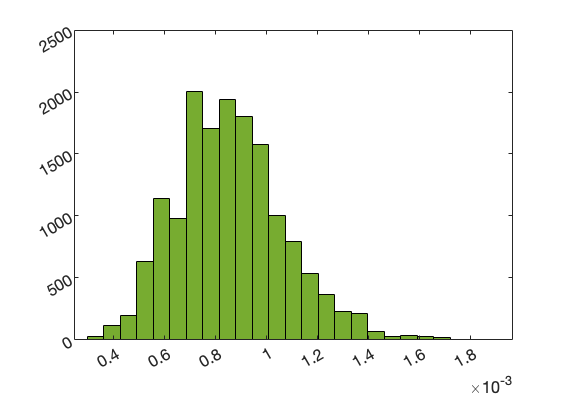}}\\
   \subfloat[Markov chain for $\beta$]{\includegraphics[width=4 cm, height = 3.2 cm]{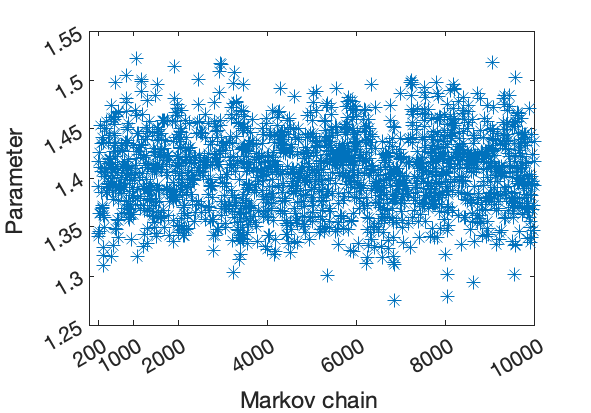}}
  \subfloat[Histogram of $\beta$]{\includegraphics[width=4 cm, height = 3.2 cm]{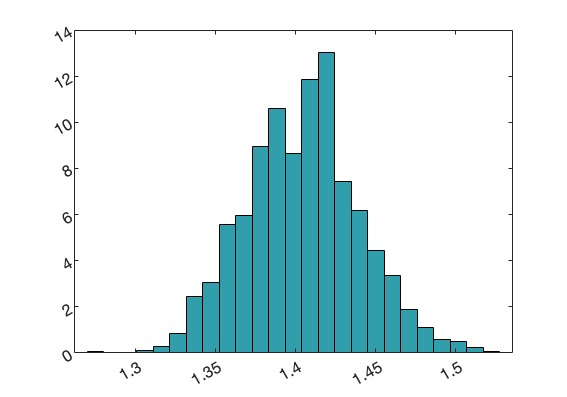}}
             \subfloat[Bivariate histogram  of $\alpha$,  $\beta$  ]{\includegraphics[width=4.2 cm, height = 3.6 cm]{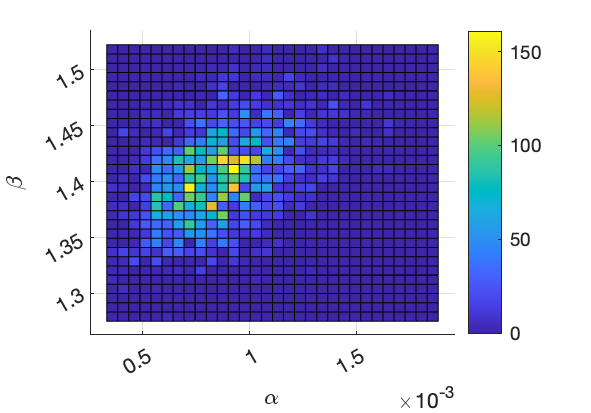}}
\caption{Diagonal solution, reference values $\alpha_{\ast}=0.0008$  and $\beta_{\ast}=1.4$, uniform prior: (a) plot of reference solution; (b) Markov chain and (c) histogram
  of the marginal posterior distribution of parameter $\alpha$; (d) Markov chain and (e) histogram of the marginal posterior distribution of parameter $\beta$;
  (f) bivariate histogram of the joint posterior distribution of $(\alpha, \beta)$.
  }
    \label{fig:D1-uniform-alpha-0.0008-beta-14-diag}
\end{figure}

\begin{figure}[H]
    \centering
      \subfloat[Reference solution]{ \includegraphics[width=0.28\linewidth]{Numerical_results/LC_example/Two_Variables/SolnD1_beta14_inL0008.png}}
    \subfloat[Markov chain for $\alpha$]{\includegraphics[width=4 cm, height = 3.2 cm]{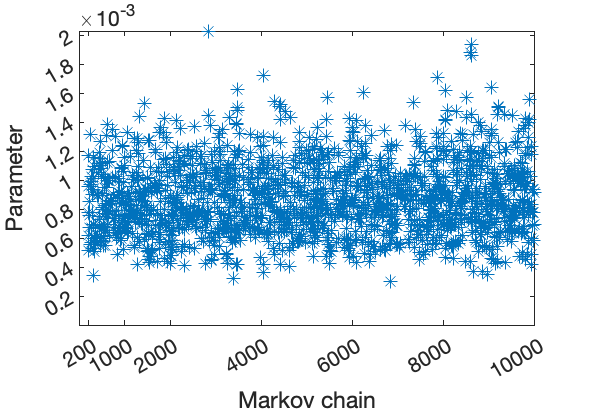}}
     \subfloat[Histogram of $\alpha$]{\includegraphics[width=4 cm, height = 3.2 cm]{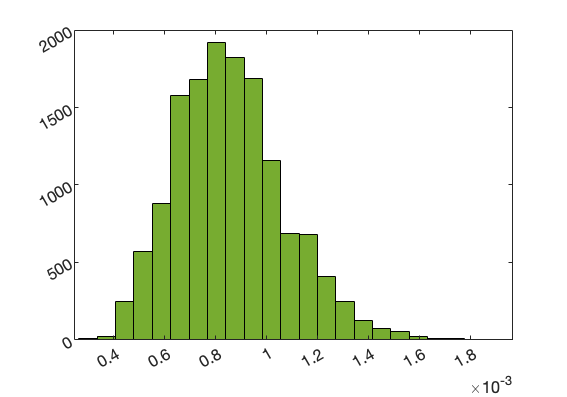}}
  \\
   \subfloat[Markov chain for $\beta$]{\includegraphics[width=4 cm, height = 3.2 cm]{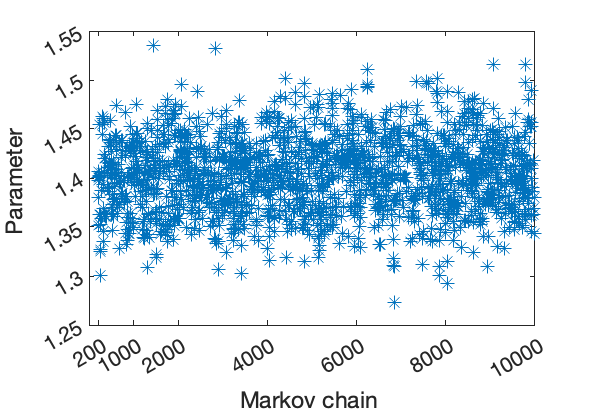}}
  \subfloat[Histogram of $\beta$]{\includegraphics[width=4 cm, height = 3.2 cm]{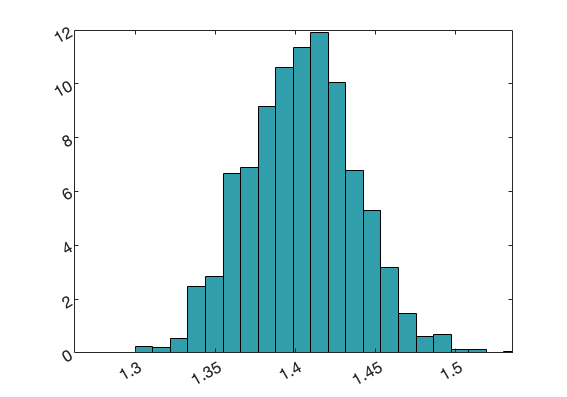}}
             \subfloat[Bivariate histogram  of $\alpha$,  $\beta$  ]{\includegraphics[width=4.2 cm, height = 3.6 cm]{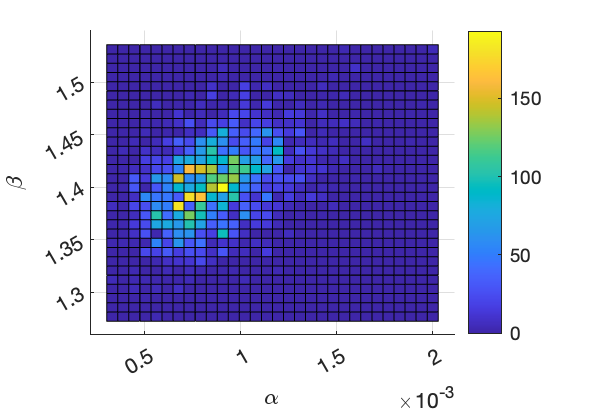}}
  \caption{Diagonal solution, reference values $\alpha_{\ast}=0.0008$  and $\beta_{\ast}=1.4$, Gaussian prior: (a) plot of reference solution; (b) Markov chain and (c) histogram
  of the marginal posterior distribution of parameter $\alpha$; (d) Markov chain and (e) histogram of the marginal posterior distribution of parameter $\beta$;
  (f) bivariate histogram of the joint posterior distribution of $(\alpha, \beta)$.
  }
    \label{fig:D1-Gaussian-alpha-0.0008-beta-14-diag}
\end{figure}

\begin{figure}[H]
    \centering
      \subfloat[Reference solution]{ \includegraphics[width=0.28\linewidth]{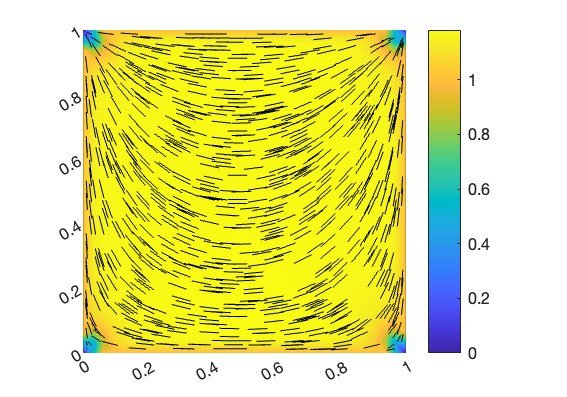}}
    \subfloat[Markov chain for $\alpha$]{\includegraphics[width=4 cm, height = 3.2 cm]{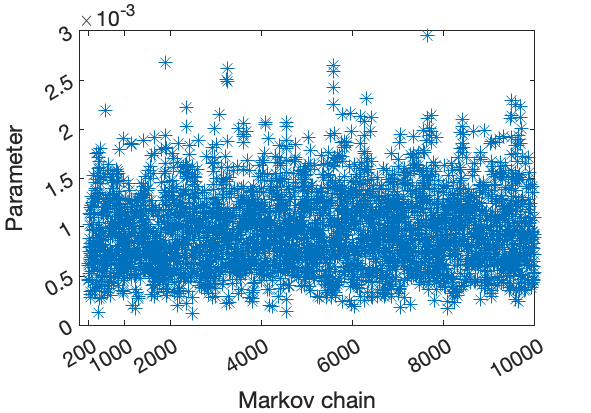}}
     \subfloat[Histogram of $\alpha$]{\includegraphics[width=4 cm, height = 3.2 cm]{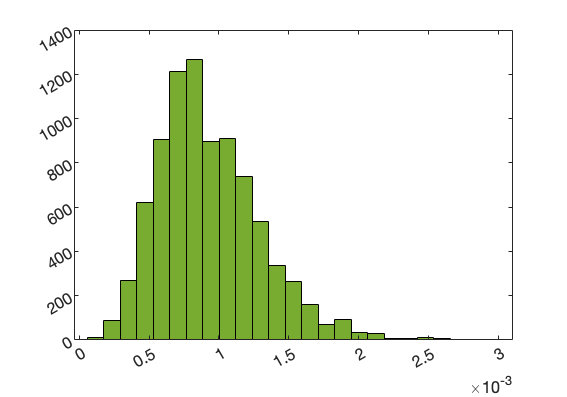}}
    \\
   \subfloat[Markov chain for $\beta$]{\includegraphics[width=4 cm, height = 3.2 cm]{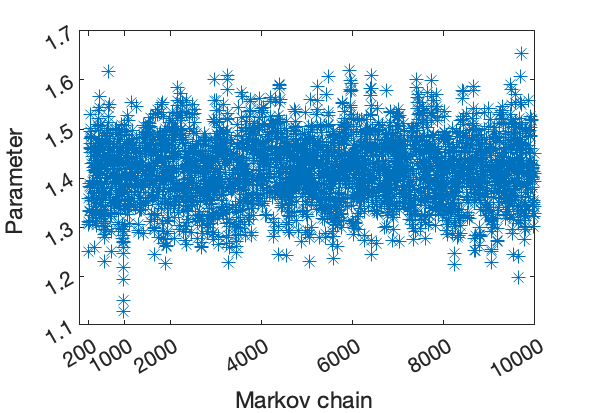}}
  \subfloat[Histogram of $\beta$]{\includegraphics[width=4.2 cm, height = 3.6 cm]{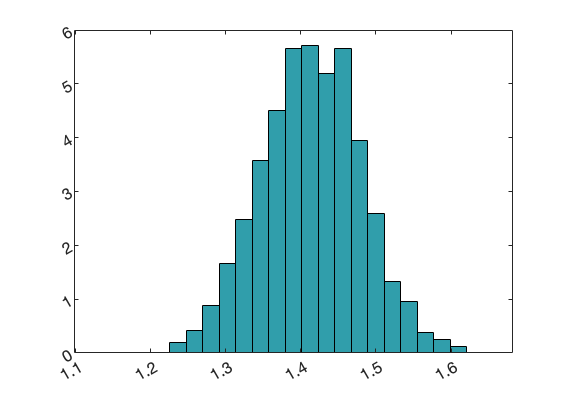}}
\subfloat[Bivariate histogram  of $\alpha$, $\beta$  ]{\includegraphics[width=4.2 cm, height = 3.6 cm]{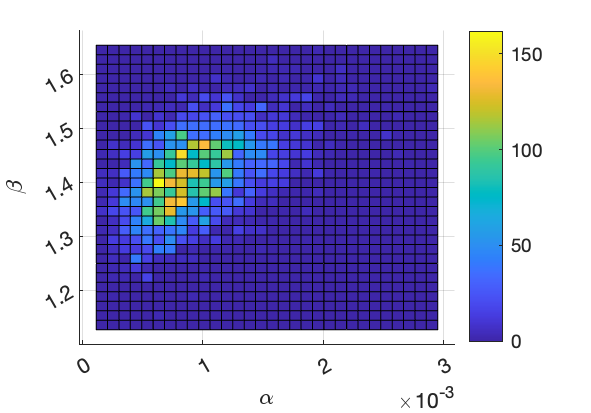}}
\caption{Rotated solution, reference values $\alpha_{\ast}=0.0008$  and $\beta_{\ast}=1.4$, uniform prior: (a) plot of reference solution; (b) Markov chain and (c) histogram
  of the marginal posterior distribution of parameter $\alpha$; (d) Markov chain and (e) histogram of the marginal posterior distribution of parameter $\beta$;
  (f) bivariate histogram of the joint posterior distribution of $(\alpha, \beta)$.
  }
    \label{fig:R4-Uniform-alpha-0.0008-beta-14-rot}
\end{figure}

\begin{figure}[H]
    \centering
      \subfloat[Reference solution ]{ \includegraphics[width=0.28\linewidth]{Numerical_results/LC_example/Two_Variables/SolnR4_beta14_inL0008.png}}
    \subfloat[Markov chain for $\alpha$]{\includegraphics[width=4 cm, height = 3.2 cm]{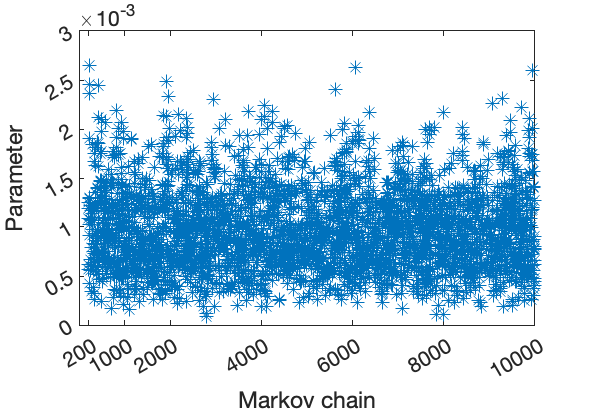}}
     \subfloat[Histogram of $\alpha$]{\includegraphics[width=4 cm, height = 3.2 cm]{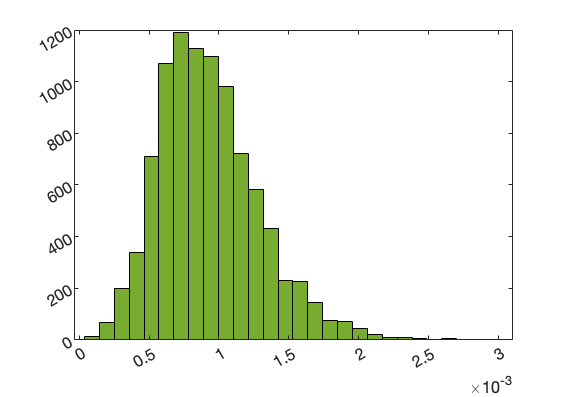}}
     \\
   \subfloat[Markov chain for $\beta$]{\includegraphics[width=4 cm, height = 3.2 cm]{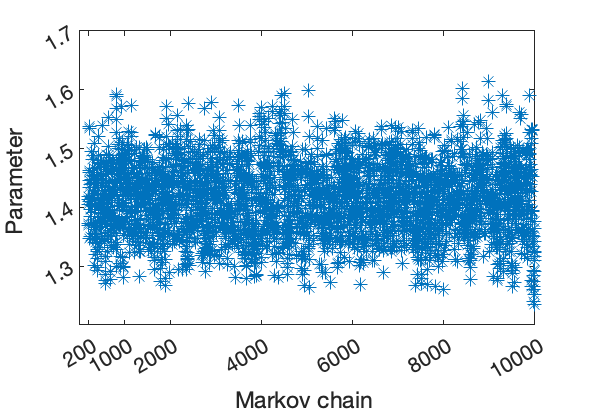}}
  \subfloat[Histogram of $\beta$]{\includegraphics[width=4 cm, height = 3.2 cm]{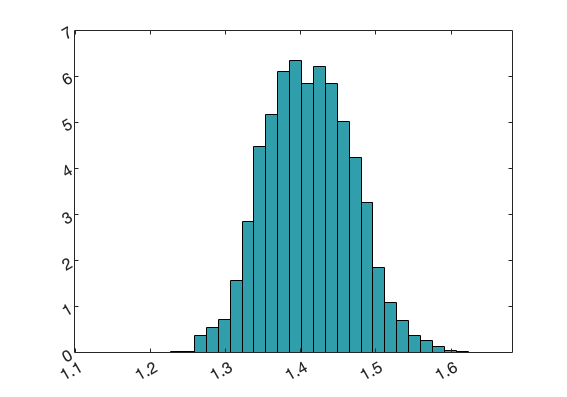}}
\subfloat[Bivariate histogram  of $\alpha$, $\beta$  ]{\includegraphics[width=4.2 cm, height = 3.6 cm]{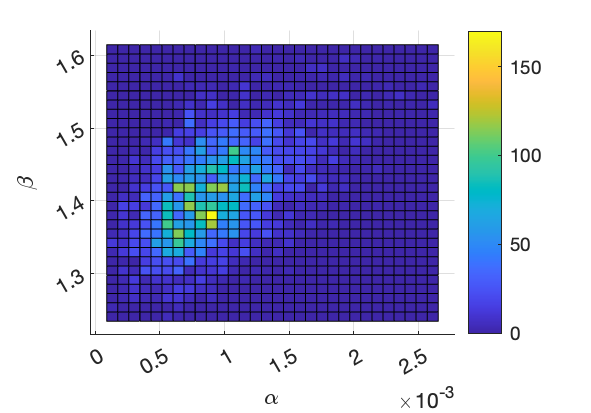}}
\caption{Rotated solution, reference values $\alpha_{\ast}=0.0008$  and $\beta_{\ast}=1.4$, Gaussian prior: (a) plot of reference solution; (b) Markov chain and (c) histogram
  of the marginal posterior distribution of parameter $\alpha$; (d) Markov chain and (e) histogram of the marginal posterior distribution of parameter $\beta$;
  (f) bivariate histogram of the joint posterior distribution of $(\alpha, \beta)$.
  }
    \label{fig:R4-Gaussian-alpha-0.0008-beta-14-rot}
\end{figure}

In Figures~\ref{fig:D1-uniform-alpha-0.0008-beta-14-diag} and \ref{fig:D1-Gaussian-alpha-0.0008-beta-14-diag}, we plot the posterior distributions for $\alpha$ and $\beta$ as described above, with a diagonal solution as $\Qvec_{\text{obs}}$. In Figures~\ref{fig:R4-Uniform-alpha-0.0008-beta-14-rot} and \ref{fig:R4-Gaussian-alpha-0.0008-beta-14-rot}, we we plot the posterior distributions for $\alpha$ and $\beta$ as described above, with a rotated solution as $\Qvec_{\text{obs}}$, when $\alpha_{\ast} = 0.0008$ and $\beta_{\ast} = 1.4$. We omit the plots for $\alpha_{\ast}=0.004, \beta_{\ast}=0.6$ for brevity.
    
\vspace*{-1.0cm}

\vspace*{-1.0cm}

\section{Examples of non-identifiability} \label{Disc}

Next, we show that the MCMC algorithm may fail to identify/reconstruct parameter values of $\alpha, \beta$, when they are outside the physical range
mentioned at the end of Subsection \ref{Inverse-problem-formulation}\ref{sec: Q to para}.
As a test example, consider the nonlinear system \eqref{pde-two-parameters} (with $\beta =1$) on a unit square domain,  $\Omega:=(0,1)^2$, with the non-homogeneous Dirichlet boundary condition
by $ \Qvec_b(x,y)=\frac{((x-a_1), (y-a_2))}{\sqrt{{(x-a_1)}^2+{(y-a_2)}^2}} .$ 
Here $a:= (a_1, a_2)$ is the centre of an interior point vortex, and we choose $a:=(1/4,3/4).$ 

\begin{figure}[H]
    \centering
     \subfloat[Reference solution]{\includegraphics[height=3 cm, width=4 cm]{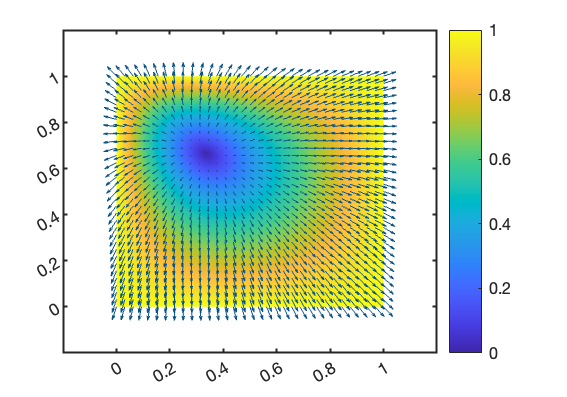}}
   \subfloat[Markov chain for $\alpha$ ]{  \includegraphics[width=4 cm, height = 3.2 cm]{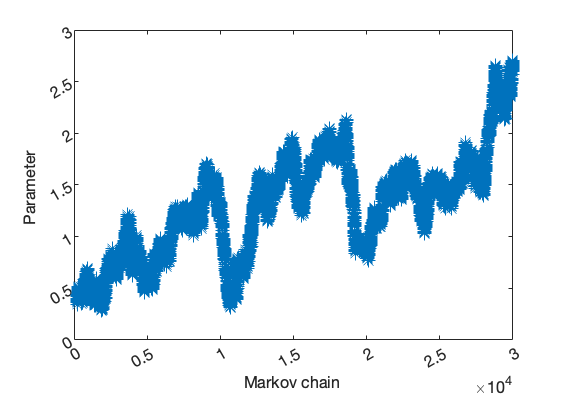}}
  \subfloat[Histogram of $\alpha$]{  \includegraphics[width=4 cm, height = 3.2 cm]{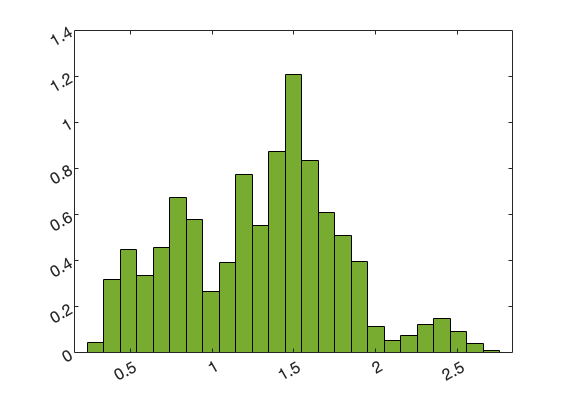}}
    \caption{(a) Computed solution, reference value {$\alpha_{\ast}=1;$} (b) Markov chain and (c) histogram of the posterior distribution  of parameter $\alpha$.}
     \label{fig:L-1-ex2}
\end{figure}

\begin{figure}[H]
    \centering
     \subfloat[Reference solution]{  \includegraphics[height=3 cm, width=4 cm]{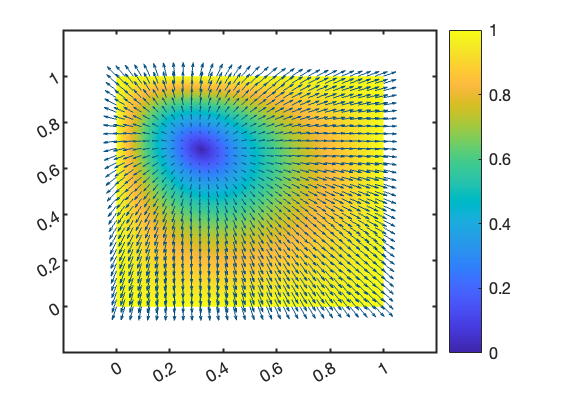}}
   \subfloat[Markov chain for $\alpha$]{   \includegraphics[width=4 cm, height = 3.2 cm]{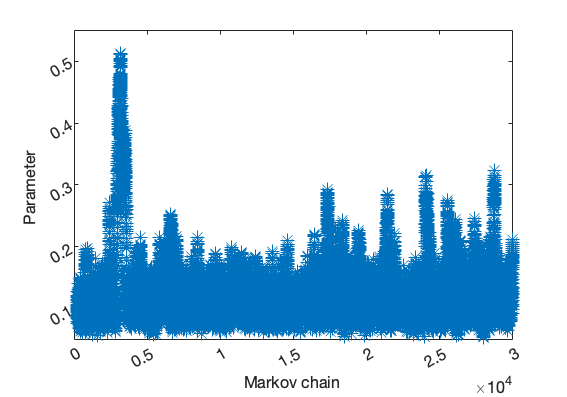}}
  \subfloat[Histogram of $\alpha$]{    \includegraphics[width=4 cm, height = 3.2 cm]{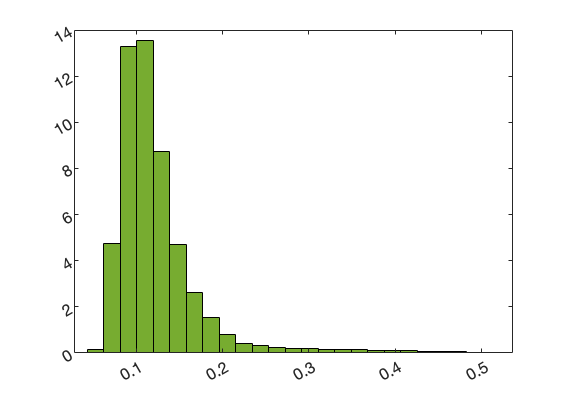}}
    \caption{
    (a) Computed solution, reference value {$\alpha_{\ast}=0.1;$} (b) Markov chain and (c) histogram of the posterior distribution  of parameter $\alpha$.
    }
    \label{fig:L-point1-ex2}
\end{figure}

\begin{figure}[H]
    \centering
     \subfloat[Reference solution]{  \includegraphics[height=3 cm, width=4 cm]{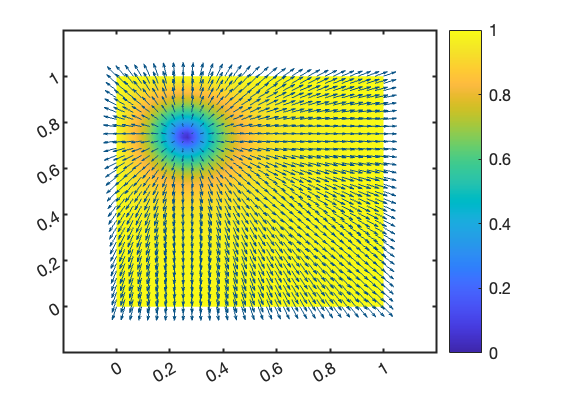}}
   \subfloat[Markov chain for $\alpha$]{   \includegraphics[width=4 cm, height = 3.2 cm]{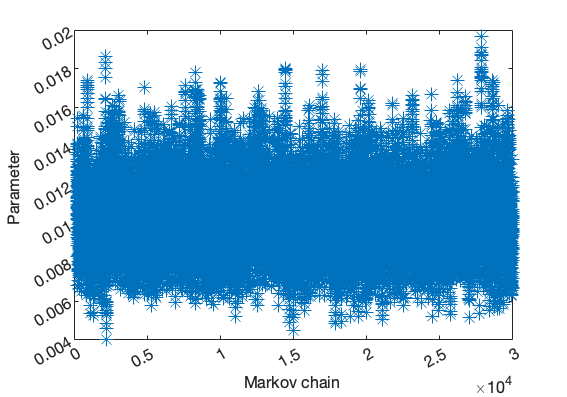}}
  \subfloat[Histogram of $\alpha$]{    \includegraphics[width=4 cm, height = 3.2 cm]{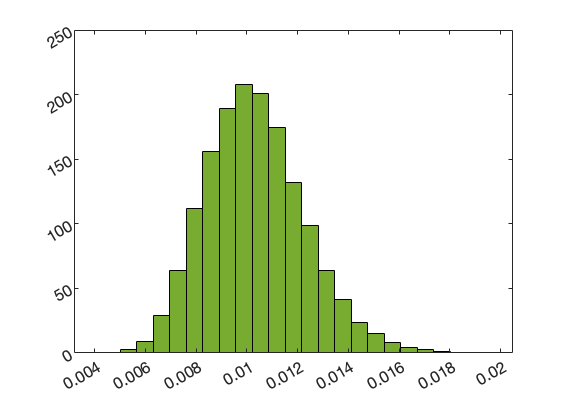}}
   \caption{
   (a) Computed solution, reference value {$\alpha_{\ast}=0.01;$} (b) Markov chain and (c) histogram of the posterior distribution  of parameter $\alpha$.
   }
    \label{fig:L-0.01-ex2}
\end{figure}

We follow the same strategy as before. We solve \eqref{pde-two-parameters} with $\alpha_{\ast} = 1, 0.1, 0.01$, $\beta_{\ast}=1$; use these solutions as the observed data $\Qvec_{\rm obs}$ in the likelihood function and then construct the posterior distribution of $\alpha_{\ast}$. Referring to Figures~\ref{fig:L-1-ex2}, \ref{fig:L-point1-ex2} and \ref{fig:L-0.01-ex2}, (a) is a plot of the numerical solution which acts as $\Qvec_{\rm obs}$.The blue arrows correspond to the vector field,
$\Qvec=(Q_{11}, Q_{12})$, and the colour-bar corresponds to $\abs{\Qvec}.$ Plot (b) is the Markov chain of length 30000 obtained by the MCMC algorithm with the improper uniform prior distribution  $\pi_{\text{prior}}(\alpha)  = \chi(\alpha)$  and the likelihood function in \eqref{likelihood-LC-example}. Plot (c) displays the corresponding histogram for the posterior distribution of $\alpha_{\ast}$.

\vspace*{-1.5cm}

\vspace*{-1.5cm}

The figures show that the Markov chain does not converge for data $\Qvec_{\text{obs}}$ constructed from $\alpha_{\ast} = 1$. It somewhat converges around $\alpha_{\ast} = 0.1$, though with a badly mixing Markov chain, and finally converges when $\alpha_{\ast} = 0.01$ (and smaller).
The profile likelihood is another method for detecting such non-identifiability  \cite{raue2013joining,siekmann2012mcmc}. The general profile likelihood can be used to detect identifiability or non-identifiability of individual components in a multi-parameter model. In our case of a scalar parameter $\alpha$, the profile likelihood is just the likelihood function as in \eqref{likelihood-LC-example}. 

\begin{figure}[h]
    \centering
      \subfloat[{$\alpha_{\ast}=1$}]{   \includegraphics[width=4 cm]{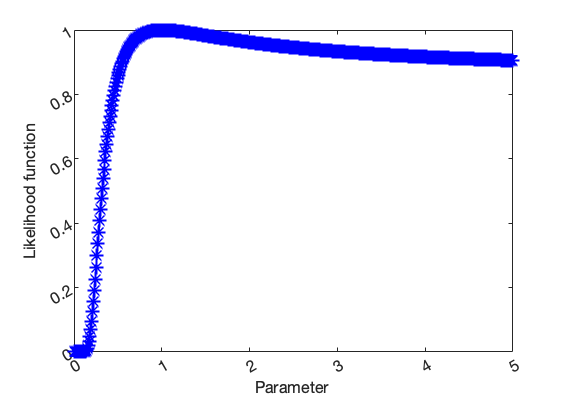}}
  \subfloat[{$\alpha_{\ast}=0.1$}]{      \includegraphics[width=4 cm]{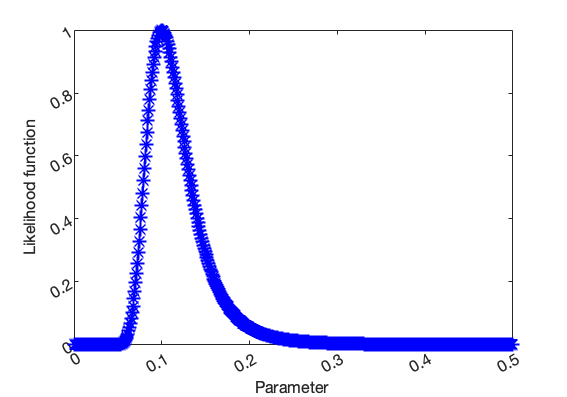}}
   \subfloat[{$\alpha_{\ast}=0.01$}]{      \includegraphics[width=4 cm]{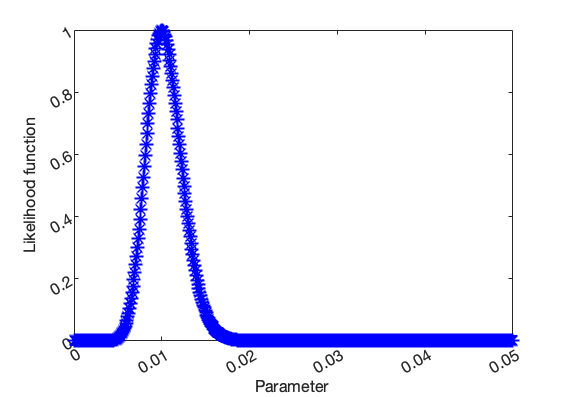}}
          \caption{Plots of the likelihood functions ${\pi_{\rm err}}(\Qvec_{\text{obs}} - \mathcal{F}(X,\alpha))$.}
    \label{fig:identifiability-ex2}
\end{figure}

If the likelihood function has a flat profile or tails out to a plateau on one or both sides, one can detect non-identifiability. This is the case in Figure \ref{fig:identifiability-ex2}(a) for $\alpha_{\ast} = 1$ and corresponds to divergence of the Markov chain as shown in Figure \ref{fig:L-1-ex2}(b). Conversely, if the likelihood function tails out to zero on both sides rather quickly, it indicates identifiability.
This is seen in Figure \ref{fig:identifiability-ex2}(c); the corresponding Markov chain, Figure \ref{fig:L-0.01-ex2}(b), shows a good mixing behaviour and converges to the respective values of $\alpha_{\ast}$, as can be seen from the resulting posterior distribution in Figure \ref{fig:L-0.01-ex2}(c).
The intermediate case is $\alpha_o = 0.1$; here the likelihood function in Figure \ref{fig:identifiability-ex2}(b) is bell-shaped, but has a rather fat tail on the right-hand side. This is reflected in the somewhat rugged Markov chain in Figure \ref{fig:L-point1-ex2}(b).

The non-identifiability for large values of $\alpha$ can be understood as follows: the system \eqref{pde-two-parameters} admits a unique solution for $\alpha$ large enough, say $\alpha \in \left(\alpha_c, \infty \right)$, or equivalently when the square edge length is small enough \cite{Karlj_Majumdar2014}. For example, the system \eqref{pde-two-parameters} has a unique solution in the range $\alpha \in \left(\alpha_c, \infty \right)$ subject to the boundary conditions in Table~\ref{table:tangential-BC}, labelled as the Well Order Reconstruction Solution (WORS). The WORS has two diagonal line defects that partition the square domain into four quadrants and the director is constant in each quadrant. In other words, different values of $\alpha$ can converge to the same director profile, making identifiability difficult for large values of $\alpha$. Further, there are questions about the validity of the LdG model on very small length scales or close to critical temperatures, as captured by large positive values of $\alpha$. The non-identifiability or divergence of the MCMC may simply be an indicator of the model failures in the physical regimes associated with large positive $\alpha$.

\section{Conclusions \label{conclusions}}

We have successfully set up an inverse UQ framework for the reduced LdG model in two dimensions, that has worked for the benchmark problem of the planar bistable nematic device reported in \cite{Tsakonas}, for physically relevant scenarios. However, this is a relatively simple benchmark problem with a square domain and regular boundary conditions. Importantly, we have tested our method with the diagonal and rotated solutions, both of which are free of interior nematic defects.

There are numerous follow-up open questions.
Will our methods work for more complex scenarios? For example, can our algorithms reconstruct LdG model inputs when $\Qvec_{\rm obs}$ has multiple interior defects? We have briefly touched on this topic in Section~\ref{Disc} but more extensive numerical experiments are required. Similar comments apply to non-convex domains or complex boundary conditions. It would be interesting to replace the Dirichlet tangent boundary conditions with surface anchoring energies (see \cite{MultistabilityApalachong}), and develop inverse UQ methods or Bayesian methods that can reconstruct the surface anchoring coefficient from $\Qvec_{\rm obs}$. We have only tested our methods with the diagonal and rotated solutions, both of which are local LdG energy minimisers. It would be interesting to see if our methods can reconstruct LdG model inputs from saddle-point solutions of the LdG Euler-Lagrange equations, many of which have symmetric constellations of interior defects \cite{hannonlinearity2021}. In future work, we plan to further develop and adapt our inverse UQ methods or Bayesian methods to more complex scenarios, for fundamental scientific curiosity and for providing new approaches for estimating material properties to the interdisciplinary liquid crystal research community.

\section*{Acknowledgements}The authors gratefully acknowledge funding from the Royal Society International Exchange Grant IES\textbackslash R2\textbackslash 242240.

\vskip2pc

\bibliographystyle{abbrv}
\bibliography{Reference_ILDG_MO_new-2.bib}

\end{document}